\documentclass{amsart}

\usepackage[left=4cm, right=4cm, top=3cm, bottom=2.6 cm]{geometry}

\usepackage{amscd,amssymb}
\usepackage{colordvi,color}
\usepackage{hyperref}

\usepackage{kpfonts}  
\usepackage[T1]{fontenc}

\usepackage{verbatim}

\usepackage{enumerate}

\usepackage{moreverb}

\newtheorem{lem}{Lemma}[section]

\newtheorem{cor}[lem]{Corollary}
\newtheorem{defi}[lem]{Definition}
\newtheorem{thm}[lem]{Theorem}
\newtheorem{prop}[lem]{Proposition}
\newtheorem{prty}[lem]{Property}
\newtheorem{fact}[lem]{Fact}
\newtheorem{ex}[lem]{Example}
\newtheorem{mainthm}{Theorem}

\makeatletter
\def\th@plain{%
  \thm@notefont{}
  \itshape 
}
\def\th@definition{%
  \thm@notefont{}
  \normalfont 
}
\def\th@procedure{
	\thm@notefont{}
	\ttfamily
}
\makeatother

\theoremstyle{definition} 

\newtheorem{rk}[lem]{Remark}

\theoremstyle{procedure} 

\newcommand{\bo}{{ 0}} 

\newcommand{\bx}{p}
 
\newcommand{\cO}{\mathcal{O}}

\newcommand{\Crit}{\Sigma} 
\newcommand{\frkm}{\mathfrak{m}}
\newcommand{\frkX}{\mathfrak{X}}
\newcommand{\mfm}{\mathfrak{m}}
\newcommand{\mult}{{\rm mult}}
\newcommand{\ord}{{\rm ord}}
\newcommand{\rdeg}{{\rm rdeg}} 
\newcommand{\Na}{\mathbb{N}}
\newcommand{\NF}{{\rm NF}}
\newcommand{\re}{\mathbb{R}}

\newcommand{\reg}{{\rm reg}}
\newcommand{\sgf}{{\sigma^*f}}

\newcommand{\sgm}{{\sigma}}

\newcommand{\dowod}{\noindent{\bf Proof:} }
\newcommand{\koniecdowodu}{\hfill$\blacksquare$\\ }
\newcommand{\Za}{\mathbb{Z}}
\newcommand{\ind}{l}

\newcommand{\mgE}{{\rho}} 
\newcommand{\mfE}{{\varphi}} 
\newcommand{\mmE}{{\mu}} 
\newcommand{\mXE}{{\nu}} 
\newcommand{\mXF}{{\varsigma}} 

\newcommand{\arc}{\gamma}
\newcommand{\arcup}{\theta}

\newcommand{\Int}{{\rm Int}}

\usepackage{enumitem} 
\setlist{nolistsep} 

\title[Relative multiplicity and degree]{Multiplicity and degree relative to a set}

\author{Vincent Grandjean}

\author{Maria Michalska}
\address{Departamento de Matem\'atica, Universidade Federal do Cear\'a
(UFC), Campus do Pici, Bloco 914, Cep. 60455-760. Fortaleza-Ce,
Brasil}
\address{Wydzia\l{} Matematyki i Informatyki, Uniwersytet \L{}\'o{}dzki, Banacha 22, 90-238 \L{}\'o{}d\'z{}, Poland}
\email{vgrandjean@mat.ufc.br}
\email{Maria.Michalska@wmii.uni.lodz.pl}

\keywords{relative degree, rectilinearization, order of vanishing, arc spaces, bifurcation values}
\subjclass[2010]{Primary 58K55, Secondary 	14P15, 32S45.}

\begin{document}

\maketitle

\begin{abstract}
The  multiplicity (resp. degree) of a function $f$ relative to a semianalytic subset~$S$ of $\re^n$  is the  greatest (resp. smallest) exponent among numbers~$j$ such that the inequality $|f(x)|\leq C\|x\|^j$  holds on $S$  near $0$ (resp. near $\infty$) for some constant~$C$. 
We show that there exists a family  of curves $\{\Gamma_d\}_{d\in \Na}$ determined only by the set such that the relative multiplicity  of any polynomial of degree $d$ is equal to its relative multiplicity with respect to $\Gamma_d$. Moreover,  a semianalytic family $(S_t)_{t\in\re^m}$ of sets  given by inequalities $f_i+t_ig_i\geq 0$ for $i=1,\dots, m$ admits a stratification of the parameter space $\re^m$  such that on each component of the top-dimensional stratum the relative multiplicity function on $\mathcal{O}_n$ does not change.  
Analogous results, assuming the data are algebraic, hold in the relative degree case.
\end{abstract}

\tableofcontents

\section*{Introduction and statement of main results}\label{sec_main_results}

The degree of real polynomial function is a rough, though simple, measure of its asymptotic behavior 
at infinity. Given  an unbounded set $S$ in $\re^n$ the  degree of a polynomial $f$ relative to the set $S$ is defined exactly in the same way\label{def_relative_degree_integer}
$$ \deg_Sf:= \min\left\{ {d}\in\Za_{\geq}: {|f(x)| \leq C\|x\|^d} {\rm \  on\  }S\setminus K  {\rm \ for\ some\  compact\ } K {\rm \ and\ constant\ }C \right\} .$$
In~\cite{NetMon} it is shown that the associated grading on the polynomial ring determines solvability of the moment problem on noncompact sets as well as existence of degree bounds in Positivstellens\"atze, compare~\cite{Schmudgen2003}. 
In general, it is difficult to decide whether the modules  $\mathcal{B}_q(S)$ of polynomials with relative degree not greater than $q$ are of finite type over the algebra $\mathcal{B}_0(S)$ of bounded polynomials, compare for instance~\cite{PlSchToric}. Therefore, in view of applications, it seems vital to provide effective methods to compute the relative degree. This is primary motivation of our results on constructive calculation of relative degree in Theorem~\ref{main_thm_TC_deg} and its invariance under change of parameters in Theorem~\ref{thm_degree_stability_20_integer}. 

The inversion of $\re^n$, taking infinity to the origin, localizes at a point the phenomena occurring at infinity.
From this point of view, the degree of a polynomial is local data at infinity about the function.
Therefore, relative degree is treated in this paper as a special case of multiplicity of a function relative to a set, which is the greatest among exponents~$w$ such that the function is bounded from above by~$\|x\|^w$ on the set, see Definition~\ref{def_mult}.  

The notion of relative multiplicity is of interest on its own. 
Explicit computation of relative multiplicity gives lower bounds on the local {\L}ojasiewicz exponent of a restriction as it is dual to relative multiplicity. Hence for example one could apply Theorem~\ref{main_thm_TC_deg} to get constructively lower bounds for the {\L}ojasiewicz exponent near a fiber in the sense of~\cite{RS}. 
Additionally, relative multiplicity is an analytic invariant. More precisely, we show that the 
set of values taken by relative multiplicity is a semi-group contained in ${1\over w}\Za_{\geq 0}$ for some positive integer $w$.
Therefore, relative multiplicity provides a grading over $\cO_n$ which is an analytic invariant for sets, see Property~\ref{prty_mult_anal_invariant}. 
In view of Theorems~\ref{thm_finite_mult} and~\ref{thm_multiplicity_20}, this provides a constructive obstruction to ambient analytic triviality in families of sets. 
This, along with several similarities  with research of~\cite{Pizza} or~\cite{CFKP}, suggests that relative multiplicity may be a bi-Lipschitz (possibly arc-analytic) invariant for sets. 

Let us present briefly two main results in the simplest form for relative degree. First, we show that computation of degree relative to a full-dimensional set reduces to considering a semi-algebraic family of curves depending on at most $n-1$ parameters. Combining this result with Proposition~\ref{prop_rectilinearization_is_polynomial}, that allows to take polynomially parametrized branches, is enough to show the following.

\begin{mainthm}[Curves testing relative degree]\label{main_thm_TC_deg}
For every semialgebraic set $S\subset \re^n$ fat at~$\infty$ there exists a sequence of semialgebraic curves
$$ \Gamma_0\subset \Gamma_1\subset \dots \subset \Gamma_d\subset \dots$$
such that for every $d$ we have
$$ \rdeg_S\equiv \rdeg_{\Gamma_d} {\it \quad on \ }\re_d[X]. $$
Moreover, the germ of $\Gamma_d$ at infinity admits at most $d\cdot { {n+d}\choose{d}}\cdot N$ branches each of which has a finite parametrization $t\to t^{-o}P(t)$, where $P$ is a univariate polynomial of degree~$k$. The positive integers $N$, $o$ and $k$ depend only on the set $S$.
\end{mainthm}
As a consequence, the computation of relative degree of a polynomial becomes in principle a symbolic operation with finite data. Contrary to general quantifier elimination methods, this discrete family of conditions has recursive description by adjoining consecutive members of a predefined analytically parametrized family, see Theorem~\ref{thm_existence_rational_TC_mult}.

We will additionally show that small changes in initial parameters of underlying set do not affect the degree function, whenever the initial parameter is generic.

\begin{mainthm}[Stability of relative degree]\label{thm_degree_stability_20_integer}
Consider two polynomial mappings $f,g:\re^n\to\re^m$. For $t\in\re^m$ put
 $$S_t:=\{f_1+t_1g_1> 0,\dots, f_m+t_mg_m> 0  \}.$$ 

There exists a nowhere dense semialgebraic subset $\mathcal{V}_{f,g}$ of the parameter space $ \re^m$ such that 
$$\deg_{S_t}=\deg_{S_s}$$
provided $s$ and $t$ lie in the same connected component of $\re^m\setminus \mathcal{V}_{f,g}$.
\end{mainthm}

This is shown by establishing a relation between the geometry of fibers at infinity on the exceptional divisor of rectilinearization of the mapping $(f,g)$ and critical values of the pullback of the mapping in Section~\ref{sec_stability}.

The paper begins with Preliminaries where we introduce relative multiplicity and its basic properties. In Section~\ref{section_resolution} we establish the tool of admissible rectilinearization and  in Section~\ref{sec_multiplicity_on_zero_divisor} continue with relevant properties of multiplicity under rectilinearization. Section~\ref{sec_TC} is devoted to showing the existence of parametric families of arcs on which every function attains its multiplicity. 
In Section~\ref{sec_stability} we study the trace of a parametrized family of sets on the exceptional divisor. With topological means we establish a relation of bifurcation values and fatness of sublevel sets in Theorem~\ref{thm_fat_sublevel_at_p}, allowing us to prove Theorem~\ref{thm_multiplicity_20}. The last section is devoted to the relative degree, where we include some consequences of interest as Theorem~\ref{thm_sublevel_set_stability} that may be of use in moment problems, we study the filtration introduced by rational relative degree on the ring of polynomials and show in Corollary~\ref{cor_rdeg_is_generically_a_finite_family} and Example~\ref{ex_there_can_be_moduli} that in multiparameter families generically one attains only a finite number of relative degree gradings but along positive codimension strata of the parameter space the relative degree may continuously change.

We aimed to make our paper self-contained and provide explicit examples. Some of our conclusions  generalize  results obtained  for example in~\cite{KMS, NetMon, MM_TC, PlSchToric}. Some theorems may be of interest on their own. For instance, there is a polynomial form of  rectilinearization  in Proposition~\ref{prop_rectilinearization_is_polynomial}; as a consequence  there exist analytic families of polynomial arcs that lift to uniform families of arcs over the resolution space whose set of end-points is dense on the zero divisor in Proposition~\ref{cor_arc_truncation_on_divisor};  one can decrease complexity  by means of essential components in Section~\ref{sec_essential_components};  and we provide a natural characterization of bifurcation values via sublevel sets in Theorem~\ref{prop_fat_sublevel_at_p}. We feel that the connection between bifurcation values and stability of relative degree deserves further investigation.

\section{Preliminaries}

\subsection{Conventions and notation} Throughout the paper  regular means real analytic.
A  regular mapping  $M \mapsto N$ is a real analytic mapping between regular manifolds $M$ and~$N$. 
Let $\cO_M$ be the sheaf of real analytic function germs on the regular manifold~$M$. 
{A} mapping is analytic on a {sub}set $S$ of $M$ if it is analytic on an open subset of $M$ 
containing $S$. Let $\cO_\bx$ be the algebra of regular function germs over $M$ at the point $\bx$ of $M$.
Let $\cO_\bx(M,N)$ be the $\cO_\bx$-module of germs at $\bx\in M$ of 
regular mappings $(M,\bx)\to N$.  
For  subsets $A ,S$ of $M$, the set of analytic arcs in $S$ with end-points in $A$ is
$$ \mathcal{L}_A(S) := \left\{\arcup\in \cO_0(\re,M) \ | \  \arcup(0)\in A {\rm\ and\ }   
\exists_{\epsilon>0}\ \theta\big((0,\epsilon]\big)\subset S\setminus A \right\}  $$

Let $X=(X_1,\dots,X_n)$ and $\re_d[X]$ be the subset of the  algebra $\re[X]$ of real polynomials in $n$ variables of polynomials of degree lower than or equal to $d$. 
A parametrization $\arc$ of an arc is Puiseux if its every coordinate is of the form $\sum_{i=0}^{\infty} a_i t^{i/q} $ 
for some some $q\in\Na$ and convergent for $0<t<< 1$. 
A parametrization~$\arc$ of an unbounded analytic arc is called Laurent-Puiseux, if its every coordinate is 
of the form $\sum_{i=-\infty}^{k} a_i t^{i/q} $, $t>> 1$ for some $k\in\Za$ and $q\in\Na$. 
Every unbounded analytic arc has such a parametrization. 
Similarly as in the case of standard Puiseux series the order $\ord \arc$ of a Laurent-Puiseux series $ \arc\in\re\{t^{1/q}, t^{-1}\}$ is the infimum of $i/q$ such that $a_i\neq 0$, and the degree $\deg \arc$ of $\arc$ is the maximum of $i$ 
such that $a_i\neq 0$.
We use convention  that $\deg 0=-\infty=-\mult 0$, where $\mult$ is the standard multiplicity at~$0$.

We will  use  Euclidean topology unless stated otherwise.  A set is fat if it is contained in 
the closure of its interior. 
 A set is fat/open/closed at $\bx\in\re^n\cup\{\infty\}$ if its germ 
at the point~$\bx$ is fat/open/closed. 
Let us denote by $\partial S$ the boundary of $S$, i.e. the set $\overline{S}\setminus \Int(S)$, and by~$\overline{A}^S$ and $ \Int_S(A)$  respectively the closure and  interior of a subset $A$ of $S$
with respect to the restricted Euclidean topology of $S$.

\subsection{Multiplicity at a point with respect to a set} Throughout the rest of this section, let $S$ be a subset of $\re^n$. We introduce the multiplicity at the origin relative to $S$ of a function germ at $0$ as a measure of the asymptotic behavior of the function at $0$ along $S$.

\begin{defi}\label{def_mult}
Let us define the multiplicity at the origin relative to the set $S$ of a regular function $f$ as
\begin{equation*}
\mult_S f := \inf_{\arc\in \mathcal{L}_0(S)} {\ord (f\circ\arc)\over \ord\|\arc\|}.
\end{equation*}
If $0\notin\overline{S}$, we put $\mult_S \equiv \infty$.
\end{defi}
Obviously, the interesting cases are when $S$ accumulates at $0$ and $S$ admits an analytic Curve Selection Lemma,in which case the following definition is equivalent.

\begin{prty}\label{prty_multiplicity_as_optimal_exponant_in_inequality}
Let $S$ be a subanalytic subset of $\re^n$. Then for any regular $f$ we have
\begin{equation*} 
\mult_S f = \sup \left\{  q\in\re \, : \,\frac{|f(x)|}{\|x\|^q} \,\mbox{\rm is bounded on }\,S\cap B(0,\epsilon) 
\,\mbox{\rm for some }\,\epsilon>0 \right\}. 
\end{equation*} 
\end{prty}

In this paper we will consider the relative multiplicity  as a function from the space of real analytic function germs at the origin
$$\mult_S : \cO_n \to\re\cup \{\infty\} $$ 
Note that the multiplicity at origin relative to a set is well-defined for arc-analytic functions germs at $0$ which includes meromorphic functions with indeterminacy at $0$. 
Note that if $S=\re^n$, then $\mult_S$ is the standard multiplicity $ \mult$ at the origin. 

\subsection{Some properties of relative multiplicity}\label{subsec_properties_of_mult}
If the ambient dimension is $n=1$, then either 
$\mult_S = \mult$ or $\mult_S $ is constant. Therefore, only the case $n\geq 2$ is of interest. 
Assume that $S$ is subanalytic.

\begin{prty}\label{prty_degree_interior_vs_closure}
We have  $\mult_S=\mult_{\overline{S}}$. Moreover,  if $S$ is fat at the origin, then $\mult_S=\mult_{\Int(S)}$.
\end{prty}

\begin{prty}
Multiplicity at the origin with respect to a nonempty germ of a subanalytic set is a valuation of $\cO_n/\mathcal{I}(S)$.
For all $f,g\in\cO_n$ we have
\begin{enumerate}
\item $\mult_S :\cO_n\to [0,\infty] $ and $\mult_Sf=\infty$ if and only if $f\equiv 0$ on $S$
\item $\mult_S (fg )= \mult_S f+\mult_Sg$
\item $\mult_S (f+g) \geq \min\{\mult_Sf, \mult_Sg \}  $
\end{enumerate}
\end{prty}

\begin{prty}\label{prop:obvious}
Let $L_\bx (S)$ be the tangent link of $S\subset\re^n$ at $\bx\in\{0, \infty\}$  
defined as
$$
L_\bx (S) :=\left\{{\bf u} \in \mathbb{S}^{n-1}: \exists_{(x_k)\subset S}\  x_k\to \bx{\ \rm and \ }\frac{x_k}{|x_k|} \to {\bf u}\right\}.
$$
If $L_0 S$ contains the germ at $0$ of an open  cone, then $\mult_S f = \mult_0 f$.
\end{prty}

\section{Rectilinearization}

Throughout the paper we use the embedded resolution of singularities, as presented below, and 
follow the language of \cite{BM_Canonical_1997,BM_IHES}. For our full purposes we work in the analytic category.

A simple normal crossing divisor (SNC) of a regular manifold $M$  is the co-support $D$ 
of a principal $\cO_M$-ideal of finite type which is locally monomial at each point of $M$ such that each 
of its irreducible components is regular.

\subsection{Admissible rectilinearization} \label{section_resolution}
Let $\mfm$ be the maximal ideal of $\cO_n$. 
Take an ideal $I$  of $\cO_n$ and $\frkX = V(I)$ its zero set germ at $0$  in $\re^n$.

Let $\sgm:(M,E\cup F,E)\to ((\re^n,0),\frkX, 0)$ be a regular proper map such that 
\begin{enumerate}
\item $M$ is a regular manifold

\item the divisor $E\cup F$ is SNC, write the zero divisor $E=H_1\cup\dots\cup H_r$ as the union of its regular irreducible components~$H_i$

\item $\sgm=\pi_w\circ\dots\circ\pi_0$, where $\pi_0$ is the blowing-up of the origin $0\in\re^n$ and 
$\pi_i$ are blowing-ups with geometrically admissible centers for $i=1,\dots, w,$

\item for any $l=1,\dots, r$ and point $\bx\in H_{i_1}\cap \dots \cap H_{i_l}\setminus \cup_{j\notin\{ i_1,\dots, i_l\}} H_j$ there exists its open neighborhood $U$ in $ M$ and a system of coordinates $(u,v) = (u_1,\dots,u_l,v_1,\dots , v_{n-l})$ at $\bx$ adapted to $E$ such that $ E \cap U = \{ u_1=\dots=u_l=0 \} $

\item the ideals $\sgm^*(\mathfrak{m})$ and $ \sgm^*(I)$ are principal and monomial i.e. there exist integer points
$\mmE, \mXE\in\Za_{\geq 0}^{r}$ such that $$\sgm^*(\mathfrak{m}) = I_{H_1}^{\mmE_1}\cdots I_{H_r}^{\mmE_r} 
\quad{\rm and}\quad  \sgm^*(I)=  I_{H_1}^{\mXE_1}\cdots I_{H_r}^{\mXE_r}\cdot I_F$$ 
where $I_{H_i}$ denotes the vanishing ideal of $H_i$ in $\cO_M$ and $I_F$ is an ideal vanishing only on $F$ 
\item the exponents $\mmE, \mXE\in \mathbb{Z}^r_\geq $ are comparable with respect to the partial ordering 
$$\mmE \succcurlyeq \mXE \iff \forall_{i=1,\dots,r} \quad \mmE_i\geq \mXE_i $$
\end{enumerate}
\label{page_order_definition_in_rectilinearization}

We will call such an embedded resolution an admissible rectilinearization of~$I$. 
If the ideal is not specified, we assume that it is the vanishing ideal of $\frkX$ and  speak of an $\frkX$-admissible rectilinearization.

\begin{fact}[\cite{Hi1, BM_Canonical_1997}]
Any germ at the origin of real analytic set admits an admissible rectilinearization.
\end{fact}

For functions $h_1,\dots, h_m\in \cO_n$ define
\begin{equation*} 
h := \prod_{i}h_i\prod_{j\neq k}(h_j-h_k).
\end{equation*}
Take an admissible rectilinearization  $\sgm$ of $\frkX := \{h=0\}$ with $E\cup F = H_1\cup \dots\cup H_r$. 
For every $i$ the ideal $\sgm^*(h_i)$ is normal crossing, i.e. 
$$\sgm^*(h_i)= I_{H_1}^{\mXE_1^i}\cdots I_{H_r}^{\mXE_r^i}$$ 
for an integer point $\mXE^i\in\Za_r$. Moreover, by~\cite{BM_IHES}  the family of exponents $\{\mu, \mXE^1,\dots,\mXE^m \}$ is totally ordered by~$\succcurlyeq$.

If it does not lead to confusion, we will often drop the indices for $H_i, \mmE_i,\mXE_i$ and write $H,\mmE,\mXE$ respectively. 

If $(E\cup F)_\reg$ is the set of regular points of $E\cup F$, 
let $$H_\reg = H\cap (E\cup F)_\reg$$ be the set of regular points of $E\cup F$ on  the component $H$. 
The nonregular points of the exceptional divisor $E\cup F$ are called corner points.

\subsection{Polynomial admissible rectilinearization}

Using the notation of the previous subsection, the following  is instrumental 
to our proof of the existence of Testing Curves in Section~\ref{sec_TC}.

\begin{fact}[\cite{BBGM}]\label{fact_special_form_BBGM}
There exists a nonempty Zariski-analytic open subset $V$ of $E$ for which the following holds: for every component $H$ of $E$ there exist nonnegative integers $d_1,\dots, d_n, q_1,\dots, q_n$ such that every point $\bx$ of $V\cap H$ admits a neighborhood $U$ in $M$ and adapted coordinates $(u,v_2,\dots, v_n)$ at 
$\bx$ such that $H\cap U= \{ u=0\}$ and the admissible rectilinearization~$\sgm$ has the following local form
\begin{eqnarray}\label{eqn_form_pi_Weierstrass_lemma}
\sgm(u,v) = \left(u^{d_1},\ u^{d_2}v_2 + u^{q_2}g_2(u), \ 
\dots,\ u^{d_n}v_n +u^{q_n}g_n(u,v_2,\dots,v_{n-1}   )   \right) ,
\end{eqnarray}
where $v=(v_2,\dots, v_n)$ and $g_j:(U,0)\to(\re,0)$ are real analytic. 
\end{fact}

For a polynomial map $u\to (\sgm_1(u),\dots, \sgm_n(u))$ its multidegree is $(\deg\sgm_1,\dots, \deg\sgm_n )$. We have the following refinement of Fact~\ref{fact_special_form_BBGM}.

\begin{prop}[Polynomial admissible rectilinearization]\label{prop_rectilinearization_is_polynomial}
Let $\sigma$ be an $\frkX$-admissible rectilinearization. There exists a nonempty Zariski-analytic open subset $V$ 
of the zero divisor~$E$ and for any component $H$ there exists an integer point $d_H\in \Za^n$ such that  any point of $H\cap V$ has an open neighborhood $U$ in $M$ in which, 
after regular change of coordinates, 
$H\cap U = \{u=0\}$,  the restriction $\sgm_{|U}$ is of the form~\eqref{eqn_form_pi_Weierstrass_lemma} and is polynomial in~$u$ with multidegree  $\deg\sgm_{|U} = d_H$.

\end{prop}

Proof follows immediately from Fact~\ref{fact_special_form_BBGM} and Lemma~\ref{claim_Weierstrass} below.

\begin{lem} 
\label{claim_Weierstrass}\label{lem_Weierstrass}

Let $\pi:U\to\re^n$ be a real analytic mapping over an open neighborhood $U$ of $0$ in $\re^{n}$ of the 
form~\eqref{eqn_form_pi_Weierstrass_lemma}. 
There exists an analytic change of coordinates such that $\pi$ remains of the form~\eqref{eqn_form_pi_Weierstrass_lemma} 
and is polynomial in $u$ with $\deg_u\pi_i=d_i$ and $\ord_u\pi_i\geq \min\{d_i,q_i \}$.

\end{lem}

\dowod
The change of coordinates can be given explicitly as the composition of 
changes~\eqref{eqn_form_of_coord_change_1} and~\eqref{eqn_form_of_coord_change_2} below.  
To prove the claim, first note that 
$$\pi_j(u,v) =  u^{d_j}v_j + u^{q_j} g_j(u, v_2,\dots, v_{j-1})  $$
depends only on the first $j$ coordinates (see \cite{BBGM}).

We use  induction. The first coordinate $\pi_1$ is  polynomial in $u$ of degree $d_1$. 
Assume that there exists an analytic change of coordinates such that coordinates $\pi_1,\dots,\pi_{j-1},$ are polynomial in $u $.

Consider two cases. First, if $d_j>q_j$, we can write 
$$ u^{q_j}g_j(u,v) = u^{q_j} a_j(u,v_2,\dots, v_{j-1}) + u^{d_j}b_j(u,v_2,\dots, v_{j-1} ) $$
with $b_j$ analytic in $u,v$ and $a_j$ polynomial in $u$ of degree $< d_j-q_j$ with coefficients analytic in $v$. 
Let $A_j:U\to U_j\subset\re^n$ be the following  analytic change of coordinates 
\begin{equation}
\label{eqn_form_of_coord_change_1}
 v_j \to v_j + b_j(u,v_2,\dots,v_{j-1}) 
\end{equation}
and identity on rest of the coordinates.
Then
$$\left( \pi_j \circ A_j^{-1}\right) (u,v) =u^{d_j}v_j+  u^{q_j} a_j(u, v_2,\dots, v_{j-1}  )  $$
is polynomial in $u$ of degree $d_j$, analytic in $v$ and  $\ord_u\pi_j\geq  q_j$.

Secondly, if $d_j\leq q_j$ we can write 
$$ \pi_j(u,v) = u^{d_j} \left( v_j + u^{q_j-d_j} g_j(u,v_2,\dots,v_{j-1})  \right) $$
In this case set the analytic change of coordinates $A_j$ as
\begin{equation}
\label{eqn_form_of_coord_change_2}
 v_j \to v_j + u^{q_j-d_j} g_j(u,v_2,\dots,v_{j-1})
\end{equation}
and identity on other coordinates. Then 
$$\left( \pi_j \circ A_j^{-1}\right) (u,v) =u^{d_j}v_j  $$
is monomial with $\deg_u\pi_j = \ord_u\pi_j = d_j$.

In both cases for $i<j$ we have $\pi_i\circ A_j^{-1} = \pi_i$  and for $i>j$ we get
$$ \left(\pi_i\circ A_j^{-1}\right) (u,v) =  u^{d_i}v_i+ u^{q_i} \tilde{g}_i(u,v_2,\dots, v_{i-1}) $$
for some analytic function $\tilde{g}_i$. 
Up to shrinking $U$, induction ends the proof.
\koniecdowodu

\section{Multiplicities and geometry on the zero divisor}\label{sec_multiplicity_on_zero_divisor}

Consider a germ at the origin of an analytic set $\frkX\subset \re^n$.
Take a proper regular mapping $\sgm: (M,E\cup F,E)\to(\re^n,\frkX,0)$ with $E\cup F$ an SNC, which is biregular between $M\setminus (E\cup F)$ and $\re^n\setminus\frkX$, $E=H_1\cup\dots\cup H_r$ with $H_i$ irreducible smooth and
$$ \sgm^*\mathfrak{m} = I_{H_1}^{\mmE_1}\cdots I_{H_r}^{\mmE_r}$$\label{page_with_max_ideal_representation_upstairs}
for the maximal ideal $\mathfrak{m}$ with $\mmE_i\in\Za_{\geq }$.

\subsection{Auxiliary notation} For any subset $A$ of $\re^n$ let 
$$A_\infty:=\overline{\sgm^{-1}\left(A \setminus  \frkX\right)}\cap E$$
be the set of accumulation points of $A$ on $E$ and
$$A_{F}:= \overline{\sgm^{-1}\left(A\cap \frkX \right)}\cap E $$
the set of accumulation points of $A\cap \frkX$ in $E$. Both $A_\infty$ and $ A_F $ are closed in
the Euclidean topology of~$M$ and $A_F\subset E\cap F$. 

Let us define \label{def_H_S}
$$ E(A):= \{H:\ H\cap A_\infty\neq \emptyset,\ H {\rm \ a\ component\ of\ } E \}.$$

We will say that the set $A$ intersects  $E$ quasi-openly if for any $H\in E(A)$
the intersection  $A_\infty\cap H$ has  nonempty interior in $H$. 
Equivalently, $A$ intersects $E$ quasi-openly if and only if the intersection of $A_\infty$ with any component of 
$E$ is either empty or Zariski-analytic dense in this component. Note that $(A\setminus \frkX)_\infty = A_\infty$, hence 
$E(A) = E(A\setminus \frkX)$ and $A$ intersects $E$ quasi-openly if and only if $A\setminus \frkX$ intersects~$E$ quasi-openly.

\begin{prty}\label{prty_fat_S_on_E}
Let $S\subset \re^n$ be semianalytic and  $\frkX$ be the {Zariski-analytic} closure of $\partial S $. 
Let $\sgm$ be a $\frkX$-admissible rectilinearization. If $S$ is fat at $0$, then $S_\infty$ is fat in $E$ and $S_F\subset S_\infty$.
\end{prty}
\dowod
Observe that $\sgm(F\setminus E)$ contains $\frkX\setminus\{0\}$. If $S$ is fat, the 
closure of the interior of $S_\infty$ is a finite union of closures of connected components of $E\setminus F$, equals $S_\infty$ and contains $S_F$. 
\koniecdowodu

\subsection{Order of vanishing and relative multiplicity}

When a coherent ideal factors as $I= I_{H}^{\mXF}\cdot I'$ with $\mXF\in\Za_{\geq 0}$ maximal, the exponent 
$ \mXF$ is the order of vanishing $\ord_H I$ of~$I$ on~$H$.

\begin{prop}\label{prop_how_to_compute_mult_via_divisor_v2}\label{prop_how_to_compute_mult_via_divisor}
Let $A\subset \re^n$, $0\in \overline{A}$. 
For any function germ $f\in \cO_n$ we have
$$ \mult_A f\geq \min\left\{ \mult_{\frkX\cap A} \, f ,\ \min_{H\in E(A)} {\ord_{H}\sgm^*f\over \ord_{H}\sgm^*\mathfrak{m}}   \right\}.  $$
 Moreover, if $\{H: A_F\cap H\neq\emptyset\}\subset E(A)$, then
$$ \mult_A f\geq \min_{H\in E(A)} {\ord_{H}\sgm^*f\over \ord_{H}\sgm^*\mathfrak{m}}  .   $$
Equality holds (in both inequalities above) if  $A$ intersects $E$ quasi-openly.
\end{prop}

\dowod  
Write $ \sgm^*\mathfrak{m} = I_{H_1}^{\mmE_1}\cdots I_{H_r}^{\mmE_r}$ and  $(\sgm^*f)= I_{H_1}^{\mfE_1}\cdots I_{H_r}^{\mfE_r}\cdot I_f$ 
with $\mfE_i\in\Za_{\geq 0}$ maximal.

Let us cover $A_\infty$ by finitely many open neighborhoods $U_p\subset M$ such that $\sgm$ has adapted coordinates on $U_p$ adapted to~$E$ at the point $p$. 
In particular $\sgm(\bigcup U_p) $ contains $ (A\setminus \frkX) \cap B(0,\epsilon)$ for some $\epsilon>0$.

Given $U_p$, up to permutation of indices, we can write $p=0\in H_1\cap\dots\cap H_k$ (note that $p$ may belong to $F\cap E$), $1\leq k\leq r$, $H_i\cap U_p = \{u_i=0 \}$ and
$$f(\sgm(u,v)) = u_1^{\mfE_1}\cdots u_k^{\mfE_k}\Psi_p(u,v) $$ with $\Psi_p$ not vanishing identically on~any component of $(E\cap U_p)\setminus F$.

Without loss of generality we can assume $${\mfE_1\over \mmE_1} = \min\left\{ {\mfE_i\over \mmE_i} : i=1,\dots,k \right\}.$$ 
Write $f\lesssim g $  if there exists a constant $C>0$ such that $f\leq C\cdot g$ and $f\cong g$ if 	$g\lesssim f\lesssim g $.
We have $(\sgm^*\mathfrak{m})_p = I_{H_1}^{\mmE_1}\cdots I_{H_k}^{\mmE_k} $ hence on~$\sgm(U_p)$, possibly shrinking $U_p$, we get
\begin{eqnarray}
 |f(x)| = |f(\sgm(u,v))| = | u_1^{\mfE_1}\cdots u_k^{\mfE_k}\Psi_p(u,v) | \lesssim | u_1^{\mfE_1}\cdots u_k^{\mfE_k} | = \nonumber \\
 = \left|   \left( u_1^{\mmE_1}\cdots u_k^{\mmE_k}  \right)^{\mfE_1\over \mmE_1}\right|  \left|u_2^{\mmE_2}\right|^{{\mfE_2\over \mmE_2} - {\mfE_1\over\mmE_1}}   \cdots \left|u_k^{\mmE_k}\right|^{{\mfE_k\over \mmE_k} - {\mfE_1\over\mmE_1}} \lesssim \nonumber\\
 \lesssim   \left|   \left( u_1^{\mmE_1}\cdots u_k^{\mmE_k}  \right)^{\mfE_1\over \mmE_1}\right|  \cong \|x\|^{\mfE_1\over \mmE_1}\nonumber
\end{eqnarray}
because ${\mfE_i\over \mmE_i} - {\mfE_1\over\mmE_1}\geq 0$ for any $i=2,\dots, k$.  Since $\sgm(U_p)$ is a finite covering of $(A\setminus \frkX)\cap B(0,\epsilon)$, we get $\mult_{A\setminus \frkX} f \geq \min\{\mfE_i/\mmE_i : H_i\cap A_\infty\neq \emptyset \}$. Recall $\mult_{A\cup B}=\min\{\mult_A,\mult_B\}$. Therefore we get the first claim. 

If $\{H: A_F\cap H\neq\emptyset\}\subset E(A)$ then we can write the above for covering $U_p$ of $A_\infty\cup A_F$ and get the second claim.

If $A_\infty\cap H_1$ has nonempty interior in $H_1$, then the minimum is attained. Indeed, if $p\in A_\infty \cap U_p \setminus V(I_f)$ and $p$ is a smooth point of $F\cup E$, then 
$$  |f(x)| \cong | u_1^{\mfE_1}| = |u_1^{\mmE_1}|^{\mfE_1\over\mmE_1} \cong \|x\|^{\mfE_1\over\mmE_1}. $$

This shows equality when $A$ intersects $E$ quasi-openly and ends the proof. \koniecdowodu

As a corollary of Proposition~\ref{prop_how_to_compute_mult_via_divisor_v2} with notation of page~\pageref{def_H_S}  
we get the following.
\begin{cor}\label{cor_equal_H_S}
Let $A,B\subset \re^n$ such that $A$ and $B$ intersect $E$ quasi-openly. Suppose that either 
$\overline{A\setminus \frkX}=\overline{A}, \overline{B\setminus \frkX}=\overline{B}$ or  $A\cap \frkX= B\cap \frkX$ 
in some neighborhood of $0$. We find  $$E(A)= E(B)\quad \Rightarrow  \quad \mult_A\equiv \mult_B. $$
\end{cor}
\dowod 
If $\overline{A\setminus \frkX}=\overline{A}$, Property~\ref{prty_degree_interior_vs_closure} gives
$$\mult_{A\setminus \frkX}=\mult_{\overline{A\setminus \frkX}} =\mult_{\overline{A}} = \mult_A. $$
Use Proposition~\ref{prop_how_to_compute_mult_via_divisor_v2} for $A\setminus \frkX$ and $ B\setminus \frkX$ to get equality. 

If $A\cap \frkX= B\cap \frkX$, use Proposition~\ref{prop_how_to_compute_mult_via_divisor_v2} for $A\setminus \frkX$
and $ B\setminus \frkX$. Since $A=(A\cap \frkX)\cup (A\setminus \frkX)$, we have 
$\mult_{A} = \min\{\mult_{A\cap \frkX}, \mult_{A\setminus \frkX} \} = \min \{ \mult_{B\cap \frkX}, \mult_{B\setminus \frkX}\} = \mult_B .$ 
\koniecdowodu

\begin{lem}[Border Lemma]\label{lem_border_v3}
Take  a smooth component $H$ of $E$. Take analytic function germs $f,g:(\re^n,0)\to \re$ and assume $\sgm^*f$ and $ \sgm^*g$  have constant multiplicities $\mfE,\mgE$, respectively, at every point of $H_\reg$. Denote $S_t=\{f+tg\geq 0\}$.

There exists an open neighborhood $V\subset M$ of $H$ and closed sets $S_+,S_-,S_=, S_o$ each of which  is closure of a union of some connected components of $V\setminus (E\cup F)$ 
such that
\begin{enumerate}
\item if $\mfE< \mgE$, then $\overline{\sgm^{-1}\left(S_t\setminus \frkX\right)\cap V }= S_+ $ for all $t>0$ 
and  
$\overline{\sgm^{-1}\left(S_t\setminus \frkX\right)\cap V }= S_- $ for all $t<0$
\item if $\mfE> \mgE$, then $\overline{\sgm^{-1}\left(S_t\setminus \frkX\right)\cap V }= S_=$ for all $t\in\re$

\item if $\partial(S_t)_\infty\cap H_\reg\neq \emptyset$, then $\mfE=\mgE$. 
\end{enumerate}

\end{lem}
\dowod 
Cover $H$ by a family $\mathcal{U}$ of open sets in $M$ such that every $U\in\mathcal{U}$ admits adapted coordinates.  Put $V=\bigcup U$. Over any such $U$ we can write  $U\cap H = \{u=0\}$ and
$$\sgm^*(f-tg) = u^\mfE \Psi_f^U - t u^\mgE \Psi_g^U $$
for some regular functions $\Psi_f^U,\Psi_g^U$, which can only vanish on $\overline{E\cup F\setminus H } $. 

Obviously 
$$\sgm^{-1}(\{f+tg=0 \}\setminus\frkX)\cap U\ =\ \{ u^\mfE \Psi_f^U - t u^\mgE \Psi_g^U =0 \}\cap U\setminus (E\cup F) .$$

Since near smooth points of $H$ the functions $\Psi_f^U,\Psi_g^U $ are units, the signs of all $\Psi_f^U,\Psi_g^U $ are constant on any connected component of $V\setminus (E\cup F)$ regardless of choice of~$U$. More precisely, for any  $p,p' $ that lie in the same connected component of $V\setminus (E\cup F)$ if $U,U'\in\mathcal{U}$ are their respective neighborhoods, then  the signs of $\Psi^U_f(p)$ and $\Psi^U_g(p)$ are equal to the signs of $\Psi^{U'}_f(p')$ and $\Psi^{U'}_g(p')$ respectively.

Assume $\mgE<\mfE$. We get locally
$$\sgm^*(f-tg) = u^\mgE \left(u^{\mfE-\mgE}\Psi_f^U - t  \Psi_g^U\right) $$
and
$$\sgm^{-1}(S_t\setminus\frkX)\cap U = \left\{u^\mgE \left(u^{\mfE-\mgE}\Psi_f^U - t  \Psi_g^U\right)\geq 0 \right\} $$
Hence the  points of $\sgm^{-1}(S_t\setminus \frkX)$ in $U\setminus (E\cup F)$ are determined by the signs of $t\Psi_g^U$ and $u^\mgE$.  Moreover, on any other component $H'$ of $E$ function $\Psi^U_g$ changes sign or keeps the sign and $u^\mgE$ does not change sign, hence it is not possible for $H'\cap \sgm^{-1}(S_t\setminus\frkX)\cap U$ to be contained in $H'$. Therefore, from constancy of sign of the units $\Psi_g^U$ on connected components, we get (1).

If $\mgE>\mfE$, we have 
$$\sgm^*(f-tg) = u^\mfE \left(\Psi_f^U - t u^{\mgE-\mfE} \Psi_g^U\right). $$
As before, we get $\sgm^{-1}(S_t\setminus \{0\})\cap U$ is either a union of components of $U\setminus (E\cup F)$ or empty but does not depend on the sign of $t$. Hence we get (2).

The last claim follows from (1) and (2). 
\koniecdowodu

\subsection{End-points of arcs on the zero divisor}\label{sec_endpoints_normal_lemma}

In this section we  state some straightforward connections between relative multiplicity and order on liftings of arcs, as well as simple but relevant consequences of polynomial form of rectilinearization concerning end-points of arcs and their liftings on the zero divisor.

\begin{prty}
If $A$ is a subset of $ \re^n\setminus \frkX$, then
 $$\mult_A f = \inf_{  }  \left\{  { \ord_0 \left((\sgm^*f)\circ \arcup\right) \over \ord_0 (\sgm^*\mathfrak{m} )\circ \arcup }   \ :\ \theta\in \mathcal{L}_{A_\infty}\left( \sgm^{-1}(A\setminus \frkX)\right)   \right\} .$$
\end{prty}

\dowod Note that there is one-to-one correspondence between the set $\mathcal{L}_0(\re^n\setminus \frkX)$ of 
nonconstant arcs  in $\re^n\setminus \frkX$ with end-point in $0$ and $\mathcal{L}_E(M\setminus (E\cup F))$, arcs in $M\setminus (E\cup F)$ with end-points in $E$. Since $\sgm$ is isomorphism outside $E\cup F$, the lift $\sgm^{-1}(\arc)$ is an arc and orders are well-defined. Moreover,  $f\circ \arc = (\sgm^*f)\circ( \sgm^{-1}(\arc) )$. 
\koniecdowodu

Write $ \sgm^*\mathfrak{m} = I_{H_1}^{\mmE_1}\cdots I_{H_r}^{\mmE_r}$ and  $(\sgm^*f)= I_{H_1}^{\mfE_1}\cdots I_{H_r}^{\mfE_r}\cdot I_{F_1}^{\mfE_1'}\cdots I_{F_l}^{\mfE_l'}\cdot I_f$ with $\mfE_i, \mfE_j'\in\Za_{\geq}$ maximal.  Denote by $\langle \cdot| \cdot\rangle$ the scalar product.

\begin{prop}\label{prop_multiplicity_via_arcs_at_corner_pints_scalar}
Let $\arc $ be an arc such that $\sigma^{-1}(\arc) \in \mathcal{L}_p (M\setminus (E\cup F)) $ with 
 $$p\in H_1\cap\dots \cap H_{s_1}\cap F_1\cap\dots\cap F_{s_2} $$ for  $s_1\geq 1, s_2\geq 0$. 
Denote $\mmE_p = (\mu_1,\dots, \mu_{s_1})$,  $\mfE_p = (\mfE_1,\dots,\mfE_{s_1})$ and $\mfE_p^F=(\mfE_1',\dots,\mfE_{s_2}')$.

There exist $a\in\Za_\geq^{s_1}, a_F\in \Za_\geq^{s_2}$ such that 
for any  $f\in \cO_n$  we have
$$\ord_0 \left(f\circ\arc \right)\geq  \langle a\ |\ \mfE_p\rangle + \langle a_F\ |\   \mfE_p^F \rangle \quad {\rm and}\quad  \ord_0 \|\arc\| = \langle a\ |\ \mmE_p\rangle  $$
with equality when $p\notin V(I_f)$.
\end{prop}

\dowod  Denote $\arcup = \sgm^{-1}\circ\arc$. Take adapted coordinates $(u,v)$ on a neighborhood $U$ of $p$ with $H_i\cap U=\{u_i=0\}$ and $F_j\cap U=\{v_j=0 \}$. We can write $\sgm^*f = u^\mfE v^{\mfE'}  \Psi$ near $p=0$. For $0<t<<1$ we have
$$ (f\circ\arc)(t) = (f\circ\sgm)\circ(\sgm^{-1}\circ \arc ) (t) = \left(u^\mfE v^{\mfE'}\cdot \Psi(u,v) \right) \circ\arcup(t) = $$
$$= \arcup_1^{\mfE_1}\cdots \arcup_{s_1}^{\mfE_{s_1}}\cdot\arcup_{s_1+1}^{\mfE_1'}\cdots\arcup_{s_1+s_2}^{\mfE_{s_2}'} (t) \cdot \Psi (\arcup(t)). $$
Function $\Psi\circ \arcup$ may vanish as $t\to 0$ only if $p\in V(I_f)$, hence we get $$\ord_0 \left(f\circ\arc \right)\  \geq\  \mfE_1\ord_0\arcup_1  +\dots + \mfE_{s_1}\ord_0\arcup_{s_1} + \mfE_{1}'\ord_0\arcup_{s_1+1}+\dots+ \mfE_{s_2}'\ord_0\arcup_{s_1+s_2}$$
with equality when $p\notin V(I_f)$. 

If we write $\sgm^*\mathfrak{m} = I_H^\mmE \cdot I $ with $\mmE$ maximal, there is a unit $\Phi$ such that $$\| x\| = \|\sgm(u,v)\| \circ \sgm^{-1}(x) = ( |u^\mmE|\cdot\Phi )\circ \sgm^{-1}(x) $$  for $x\in \sgm(U\setminus (E\cup F))$. Hence $\|\arc\| = |\arcup_1^{\mmE_1}\cdots \arcup_{s_1}^{\mmE_{s_1}}| \cdot \Phi(\arcup) $ and 
$$\ord_0 \|\arc\| = \mmE_1\ord_0\arcup_1  +\dots + \mmE_{s_1}\ord_0\arcup_{s_1}.$$ 
This ends the proof.
\koniecdowodu

As a consequence we get

\begin{prty}\label{prty_how_to_get_mult_H_by_arcs}
Let $\arc:(0,1] \to\re^n $ be an arc such that $ \arc((0,1])_\infty =\{p\} \subset H_\reg $ for some component of $E$. 
For any regular $f$ we have
$$ {\ord_0 f\circ\arc \over \ord_0\|\arc\|} \geq  { \ord_H \sgm^*f\over \ord_H \sgm^*\mathfrak{m}   }   $$
with equality when $p\notin V(I_f)$.
\end{prty}

\dowod By Proposition~\ref{prop_multiplicity_via_arcs_at_corner_pints_scalar} there exists $a\in\Za_\geq$ such that for any regular $f$ we get
 $$\ord_0 f\circ\arc \geq a\cdot\ord_H \sgm^*f   \quad{\rm and}\quad \ord_0 \|\arc\| = a\cdot\ord_H \sgm^*\mathfrak{m} $$
 with equality when $p\notin V(I_f)$.
\koniecdowodu

\begin{lem}[Normal arc lemma]\label{lemma_normal_are_testing_curves}
Take a component $H$ of the zero divisor. Consider any open set $U \subset M $ such that $U\cap H=U\cap H_\reg$ is nonempty and $U$ admits adapted coordinates $(u,v)\in\re\times\re^{n-1}$ in which $H\cap U = \{u=0\}$. Without loss of generality, we may assume $U=(-\epsilon,\epsilon)\times \tilde{V}$, where $\tilde{V}$ is open neighborhood of $0$ in $\re^{n-1}$.

Consider the family
$$ \arcup^H_v =\{ (t,v)\ :\ t\in(0,\epsilon) \}.$$
We have 
$$ \ord_H \sgm^*\mathfrak{m} = \ord_0 (\sgm \circ\arcup_v^H). $$

Moreover, for any regular function $f\in\cO_n$ there exists a Zariski closed subset~$\Xi_f$ of $\tilde{V}$ such that for every $v\in \tilde{V}\setminus \Xi_f$ we have
$$ \ord_H \sgm^*f = \ord_0 f(\sgm \circ\arcup_v^H) .$$

\end{lem}
\dowod
Take $U$ satisfying the assumptions. Note that the order $\ord_{0} (\arcup_v^H)_1$ of the first coordinate of $\arcup_v^H$ is $1$. 
 As usual we can write $\sgm^*(f) = I_H^{\ord_H f}\cdot I_f $ with $I_f$ an ideal vanishing on at most a nowhere dense subset of~$H$. Hence function $\sgm^*f$ has constant multiplicity on $\tilde{V}\subset H_\reg$ outside a  nowhere dense regular subset $\Xi_f:=V(I_f)\cap \tilde{V}$. 
Therefore,  for all $v\in \tilde{V}\setminus \Xi_f$ by Property~\ref{prty_how_to_get_mult_H_by_arcs} we have
$$ \ord_H \sgm^*f =   \ord_0 f(\sgm (\arcup_v^H)) \quad {\rm and }\quad \ord_H \sgm^*\mathfrak{m} = \ord_0 \sgm (\arcup_v^H),  $$
which ends the proof.  \koniecdowodu

Consider the space\label{page_defi_Ld}
$$ \mathcal{L}^d := \left\{  \sum_{i=1}^d \overrightarrow{a}_i t^i :  \overrightarrow{a}_i\in\re^n  , t\in(0,1] \right\} $$
of polynomial mappings $(\re,0)\to (\re^n,0)$  of degree at most~$d$. It is the truncation up to degree $d$ of the space  $\mathcal{L}_0(\re^n)$ of analytic arcs in $\re^n$ with end-points at $0$.

\begin{lem}\label{lem_normal_lem_is_polynomial}
Assume that rectilinearization $\sgm$  is polynomial of the form~\eqref{eqn_form_pi_Weierstrass_lemma}. Then, under notation of Lemma~\ref{lemma_normal_are_testing_curves}, there exists $\Delta_H$ such that 
$$ \sgm(\arcup_v^H) \in \mathcal{L}^{\Delta_H}$$
for every $v\in \tilde{V}$.
\end{lem}

\dowod Due to Lemma~\ref{lem_Weierstrass} we may assume $\sgm$ is polynomial in $u$ with degree $\deg_u\sgm=\max_i d_i $. The family of parametrised arcs $\Gamma_v:=\sgm\circ\arcup_v^H$ is a subset of~$\mathcal{L}^{\max_i d_i}$. Hence $\Delta_H:=\max_i d_i$ satisfies the claim for the open subset $\tilde{V}$ of $H$. \koniecdowodu

As a consequence we get
\begin{prop}[Arc truncations]\label{cor_arc_truncation_on_divisor}
For every irreducible component $H$ of $E$ there exists a positive integer $\Delta_H$ such that for every $d\geq \Delta_H$ the set 
$$ \left\{  \lim_{t\to 0^+} \left(\sgm^{-1}\circ\Gamma\right)(t) : \Gamma\in\mathcal{L}^{d}       \right\} $$
is a dense open subset of $H$.
\end{prop}
\dowod  
Under notation of Lemma~\ref{lemma_normal_are_testing_curves} obviously  the set  end-points of $\arcup_v^H$ equals~$\tilde{V}$. Note that for arcs with finite expansions we can set the domain of their parametrization to be fixed. 
By Proposition~\ref{prop_rectilinearization_is_polynomial} the rectilinearization $\sgm$ 
may be written in the special form polynomial in $u$ in some neighborhood of any $p$ in an open dense subset $V$ of $E$ and $\Delta_H$ is independent of the choice of $U$ admitting adapted coordinates. Hence Lemma~\ref{lem_normal_lem_is_polynomial} holds on open sets $\tilde{V}$ whose union is equal to $V$.  This ends the proof. \koniecdowodu

Note that properties in this section hold under milder assumptions on the rectilinearization~$\sgm$ than the full list of Section~\ref{section_resolution}.

\subsection{Essential components of the zero divisor}\label{sec_essential_components}

In order to compute $\mult_S$ we would like to make the following expected observation: 
it may happen that to that  purpose some blowings-up are unnecessary. We will see that only certain 
components of the zero divisor, which we call essential, are pertinent to  {compute the} relative multiplicity.

Let $\sgm$ be an admissible rectilinearization of a  {real analytic set germ $\frkX$ at the origin of $\re^n$}.
We keep up with the notations of the previous sections and subsections, except that for the purpose 
of this subsection $H_1,\ldots, H_q$  {denote} the components of the divisor $E\cup F$.  {We also} assume 
that the regular hypersurfaces $H_j$ are indexed accordingly to their year of birth in the resolution process. 
(This means we assume that we blow-up connected regular center rather than the whole geometrically admissible 
center given by the algorithm which may have finitely many connected components.)

Let $C_0 = \bo, C_1, \ldots,C_\ind$ be the successive regular centers of blowing-up with $\ind \leq q-1$. 
If $C_i$ has several irreducible components, blowing-up any two such irreducible components 
is,  {up to an isomorphism}, a commutative operation, since they are disjoint. 
Thus we assume that each center blown-up was irreducible, so that $\ind=q-1$.
Let $H_i$ be the exceptional hypersurface created by blowing-up the center $C_{i-1}$, where $i=1,\ldots,q$. 
Let 
$$
E^i\cup F_i := \cup_{j=1}^i E_j \ \mbox{ for }\ E^i := \sgm_i^{-1}(0)
$$
where $\sgm_i:M_i\mapsto\re^n$ is the composition of the successive blowings-up of centers $C_0,\ldots,C_{i-1}$.


\begin{defi}[Essential components]\label{def:essential} 
A component $H_k$ of $E\cup F$ is called  $S$-inessential  if 
\begin{itemize}
\item either $H_k$ is not contracted onto $\bo$ by $\sgm$,
\item or $H_k$ is contracted to $\bo$ by $\sgm$ and the center $C_{k-1}$ is contained exactly in the components
$H_{j_1},\ldots,H_{j_m}$ of $E^{k-1}$, for some $m\geq 1$, such that 
for each connected component $\mathbf{S}$ of the germ $S\setminus \{\bo\}$ and each $j = j_1,\ldots,j_m$, 
the set $$\overline{\sgm_{k-1}^{-1}(\mathbf{S})\setminus E^{k-1}}\cap H_j$$ is either 
Zariski-analytic dense in $H_j$ or empty.
 
\end{itemize}
A component $H$ of $E$ is  $S$-essential  if it is not $S$-inessential and $S_\infty$ is Zariski-analytic dense in~$H$.
\end{defi}

We obtain a version of 
Proposition \ref{prop_how_to_compute_mult_via_divisor} as follows.

\begin{cor}\label{cor_essential_mult} Let $S$ be a  {semianalytic} subset of $\re^n$, 
 {fat at the origin}, and $\sgm$  {be} an 
admissible rectilinearization of~$\frkX$~, the  Zariski-analytic closure of  {$\partial S$}. 
Assume $\overline{S} = \overline{S\setminus \frkX} $. 

For any function germ $f\in \cO_n$ we have
\begin{equation*}
 {\mult_S} f := \min\left\{\frac{\ord_H \sgf}{\ord_H\sgm^*\frkm}\ :\ H \;
\mbox{\rm  is $S$-essential}
\right\}.
\end{equation*}
\end{cor}

Proof of Corollary~\ref{cor_essential_mult} is a consequence of Lemma~\ref{lem:key-observation} below.

\begin{lem}\label{lem:key-observation}
Let $H$ be a component of $E^i$. 
Let $H=D_1$, $D_2,\ldots,D_m,$ be all the components of $E^k$ containing 
the center $C_k$ for $k\geq i$. 
Then for any function germ $f$ of~$\cO_n$ we have
$$
\frac{\ord_{D_{k+1}} \sgm_{k+1}^* f}{\ord_{D_{k+1}} \sgm_{k+1}^*\frkm} \; \geq \;
\min_{i=1,\ldots,m}\frac{\ord_{D_i} \sgm_k^* f}{\ord_{D_i} \sgm_k^*\frkm}
$$
\end{lem}

\dowod
Fix $f$. 
Thus we can write 
\begin{eqnarray*}
\sgm_k^* (f)  =  I_{D_1}^{\mfE_1} \cdots I_{D_m}^{\mfE_m} \cdot I_f, \quad
\sgm_k^* (\frkm)  =  I_{D_1}^{\mmE_1} \cdots I_{D_m}^{\mmE_m} \cdot I_0
\end{eqnarray*}
for non-negative integer numbers $\mmE_1,\ldots,\mmE_m,\mfE_1,\ldots,\mfE_m$ and for ideals 
$I_f,I_0$ with multiplicity zero along each $D_i$.
If $D_i$ is not a component of $E^k$ then $\mmE_i = 0$.  {Since the 
co-support} of $\sgm_k^*(\frkm)$ is $E^k$, the ideal $I_0$ is a product of powers of the other  {ideals} of 
components of $E^k$ which do not contain $C_k$.
Let $\beta: (M_{k+1},D_{new}) \mapsto (M_k,C_k)$ be the blowing-up with center $C_k$. 
 {Since $\sgm^*_k(\frkm)$ is principal and monomial in $E^k$,
the multiplicity of $I_0$ along~$C_k$ is $0$. Let $\mfE$ be the multiplicity of $I_f$ along $C_k$.}
Since $\sgm_{k+1} = \beta\circ\sgm_k$, we get (using the same notation to write the strict transform by $\beta$
of each divisor $D_i$)
\begin{eqnarray*}
\sgm_{k+1}^* (f) & = & (I_{D_1}^{\mfE_1} \cdots I_{D_m}^{\mfE_m}) \cdot I_{D_{new}}^{\mfE + \sum_i\mfE_i}\cdot I_f'\\
\sgm_{k+1}^* \frkm & = & (I_{D_1}^{\mmE_1} \cdots I_{D_m}^{\mmE_m}) \cdot I_{D_{new}}^{\sum_i \mmE_i} \cdot I_0'
\end{eqnarray*}
so that
\begin{equation*} 
\frac{\ord_{D_{new}} \sgm_{k+1}^* f}{\ord_{D^{k+1}} \sgm_{k+1}^*\frkm} = \frac{\mfE + \sum_i\mfE_i}{\sum_i \mmE_i} 
\geq \frac{\sum_i\mfE_i}{\sum_i \mmE_i} \geq 
\min_{i=1,\ldots,m}\frac{\mfE_i}{\mmE_i} = 
\min_{i=1,\ldots,m}\frac{\ord_{D_i} \sgm_k^* f}{\ord_{D_i} \sgm_k^*\frkm}
\end{equation*}
This ends the proof.
\koniecdowodu

We also get the following criterion of comparability of the relative  {multiplicity} functions.

\begin{prop}\label{prop_mult_equal_wrt_essential_divisor}
Let $S_1$ and $S_2$ be semianalytic subsets of $\re^n$ both fat at the origin.
Let $\frkX$ be the Zariski-analytic closure of $\partial S_1 \cup \partial S_2$ and $\sgm$ be a
 {admissible rectilinearization $\frkX$}. 
The following are equivalent:
\begin{enumerate}
\item $ {\mult_{S_1}}\equiv  {\mult_{S_2}}$,
\item every component $H$ of $E$ is $S_1$-essential if and only if it is $S_2$-essential.
\end{enumerate}
\end{prop}

\dowod
It is a consequence of Property \ref{prty_fat_S_on_E}, 
Definition \ref{def:essential} and Corollary \ref{cor_essential_mult}.
\koniecdowodu

\subsection{Further properties of  relative multiplicity}

The local algebra of the regular function germs at the origin of a subset $S$ of $\re^n$, is the algebra 
of the germs at the origin of the regular function restricted to $S$, and 
coincides with $\mathcal{O}_{\overline{S}^Z}$ the local algebra of regular germs of the Zariski-analytic closure of~$S$. 
For brevity we denote it $\mathcal{O}_S$. {If} $S$ is fat, then $\mathcal{O}_S=\mathcal{O}_n$.

\begin{prty}
Multiplicity  relative to a semianalytic set $S$ is a discrete valuation 
on its local algebra of regular germs. 
\end{prty}

\dowod
Take a semianalytic set $S$ and the admissible rectilinearization $\sgm$ of the Zariski closure of $ \partial S $. 
Under notation as on page~\pageref{page_order_definition_in_rectilinearization}, since the number of components of the zero divisor $E$ is finite, by Proposition~\ref{prop_how_to_compute_mult_via_divisor} there exists $w\in\Na$ such that $\mult_S f\in {1/w}\Za_{\geq }$ for all $f\in\mathcal{O}_S$. This $w$ can be taken as greatest common denominator of multiplicities $\mmE_1, \dots, \mmE_r$ of the principal ideal $\sgm^*\mathfrak{m}$ on essential components of the zero divisor.
\koniecdowodu

Therefore, relative multiplicity introduces a grading 
$$ \mathcal{O}_S = \bigoplus_{k\in\Za_{\geq}} \mathcal{M}_{{k\over q(S)}}(S), $$ 
where  $ \mathcal{M}_q(S)  = \{ f\in\mathcal{O}_S : \mult_S f = q \}$ and $q(S)$ is the smallest positive integer that admits this grading.  
Moreover, relative multiplicity is an analytic invariant in the following sense.

\begin{prty}\label{prty_mult_anal_invariant}
Let $S$ and $T$ be semianalytic subsets of $\re^n$. If there exists a biregular map 
germ $ \phi : (\re^n,S,0)\to (\re^n,T,0)$, then $q(S)=q(T)$ and there exists an isomorphism of local algebras between 
$\mathcal{O}_S$ and $\mathcal{O}_T$  which preserves the gradings by {the} respective relative multiplicities.
\end{prty}
\dowod
The mapping $\phi^*:\mathcal{O}_T \to \mathcal{O}_S$ is the isomorphism. 
Obviously, $\phi^*$ is surjective. Moreover, a bianalytic $\phi$ is in particular bi-Lipschitz in a neighborhood of the origin. 
Now the proof follows easily from the fact that due to Property~\ref{prty_multiplicity_as_optimal_exponant_in_inequality} multiplicity relative to $S$ of function $f$ can be interpreted as the optimal $q\in\mathbb{Q}$ such that $|f(x)| \precsim \|x\|^q $ on $S$. \koniecdowodu

\section{Testing curves}\label{sec_TC}

This section aims at exhibiting in Theorem~\ref{thm_existence_rational_TC_mult}  a very small subset, namely a  family of polynomially parameterized arcs, so-called Testing Curves, along which the relative multiplicity of any function germ along a semi-analytic subset $S$ fat at $0$ is reached. Testing curves provide an efficient symbolic calculation  method of computing the multiplicity with respect to~$S$. In Theorem~\ref{thm_finite_mult} we will show that to study relative multiplicities of polynomials it suffices to take a countable subfamily of Testing Curves.

\subsection{Definition and existence} Let us introduce the notion of testing curves. By a curve we mean a purely one-dimensional set.

\begin{defi}[Testing Curves for relative multiplicity]\label{def_TC_mult}
For  a set  $S$ in $\re^n$ let $(\Gamma_c)_{c\in U}$, where $U\subset\re^m$, be a family of curves such that
\begin{enumerate}
\item there are $l$ connected components $\arc_1^c,\dots,\arc_l^c$ of $\Gamma_c$ for every $c$, and each component is parametrized as an analytic arc
\item  the family $(\arc_i^c)_{c\in U}$ is parametrized analytically by $c$ for each $i\leq l$
\item for every real analytic $f:(\re^n,0)\to(\re,0)$ there exists a nowhere dense Zariski closed set $\Xi_f\subset U$ such that for $c\notin \Xi_f$ we have $$\mult_S f = \mult_{\Gamma_c}f  $$ 
\end{enumerate}
Whenever $m<n$, we will call such a family $(\Gamma_c)$ the family of curves testing multiplicity with respect to $S$. 
\end{defi}
Note that if Testing Curves exist, for generic $c$ one calculates $\mult_S$ symbolically by
$$\mult_{\Gamma_c}f = \min_{i} {\ord_0 (f\circ\arc_i^c) \over \ord_0 \|\arc_i^c\| } . $$ 
Moreover, we do not ask that testing curves lie in $S$.

\begin{thm}[Existence of  curves testing multiplicity]\label{thm_testing_curves_mult_v2}
Let $S$ be a semianalytic subset of $\re^n$ fat at~$0$. There exists a family of testing curves for multiplicity at $0$ with respect to~$S$.
\end{thm}

\dowod Let $\frkX$ be the germ of the Zariski closure of $\partial S$. Consider a $\frkX$-admissible rectilinearization~$\sgm$. 
In particular $S_\infty$ is fat in~$E$ and $S_F\subset S_\infty$, see Property~\ref{prty_fat_S_on_E}. 
Under notation of Lemma~\ref{lemma_normal_are_testing_curves} there are $l$ components $H_i$ of the zero divisor $E$ such that $H_i\cap S_\infty\neq \emptyset$. Let $\arc_v^i:=\sgm(\arcup_v^{H_i})$, $v\in \tilde{V}_i$, $t\in(0,\epsilon_i)$.  
Put $$U=\bigcap_{i\ :\ H_i\cap S_\infty\neq\emptyset} \tilde{V}_i $$ 
and $\Gamma_c$ as the union of images of these arcs, i.e.
$$\Gamma_c = \bigcup_{i\ :\ H_i\cap S_\infty\neq\emptyset} \arc_c^i, \quad c\in U.$$  
Moreover, set $S$ intersects $E$ quasi-openly and  $H_i\cap S_\infty\neq \emptyset$ if and only if $i\in E(S)$. 
Hence the family $\Gamma_c$ is a family of testing curves for multiplicity relative to~$S$ by Property~\ref{prty_fat_S_on_E}, Lemma~\ref{lemma_normal_are_testing_curves} and Proposition~\ref{prop_how_to_compute_mult_via_divisor}.  \koniecdowodu

\begin{rk}\label{rk_number_of_branches_TC}
Let $\frkX= \overline{\partial S}^{\rm Zar}$ and choose $\sigma:(M,E\cup F, E)\to (\re^n, \frkX, 0)$ an $\frkX$-admissible rectilinearization with $H_1,\dots, H_r$ being all components of $E$. The number $l$ of branches of testing curves $\Gamma_c$ satisfies $$l\ \leq \ \#\{i \ :\ H_i\ {\rm is\ } S-{\rm essential} \}\ \leq\  \#\{ i\ :\ H_i\cap S_\infty\neq\emptyset\} = \#E(S).$$
\end{rk}

\subsection{Polynomially parametrized testing curves}

We show that one can always choose testing curves to have branches parametrized by polynomials whose  degree depends on $S$ and choice of rectilinearization. In particular Theorem~\ref{thm_existence_rational_TC_mult} below implies that for parametrization of testing curves it suffices to look in a finite dimensional vector space~$\mathcal{L}^d$ of truncated arcs, see page~\pageref{page_defi_Ld}, and the testing curves are semialgebraic.

\begin{thm}[Existence of polynomial Testing Curves]\label{thm_existence_rational_TC_mult}
For every semianalytic subset~$S$ of $ \re^n$ fat at  $0$ there exists an analytically parametrized family of semialgebraic testing curves for multiplicity with respect to the set $S$.

Precisely, for the semianalytic fat set $S$ there exist positive integers $\ord \Gamma, \deg \Gamma$ and testing curves satisfying Definition~\ref{def_TC_mult} with point (1)  replaced by
\begin{itemize}
\item[(1')] there are $l$ connected components $\arc_1^c,\dots,\arc_l^c$ of $\Gamma_c$ for every $c\in U$ and each component is  parametrized by a finite series
$$ \arc_j^c(t) = \sum_{i=\ord\Gamma}^{\deg\Gamma} \overrightarrow{\varrho}_{j,i}(c) t^i $$
where $\overrightarrow{\varrho}_{j,i}:U\to\re^n$ are real analytic. 
\end{itemize}
\end{thm}
\dowod Proof is the same as of Theorem~\ref{thm_testing_curves_mult_v2} but we use the special admissible rectilinearization of Fact~\ref{fact_special_form_BBGM} and apply Lemma~\ref{lem_normal_lem_is_polynomial}.\koniecdowodu

\begin{rk}\label{rk_ord_deg_Gamma_TC}
In Theorem~\ref{thm_existence_rational_TC_mult}, under notation of Fact~\ref{fact_special_form_BBGM}, by Lemma~\ref{lem_Weierstrass} we get 
$$ \ord\Gamma \geq \min_{i=1,\dots, n}\{d_i,q_i\} \quad {\rm and} \quad \deg\Gamma = \max_{i=1,\dots, n} d_i .$$
\end{rk}

\begin{ex}\label{ex_3dim}
Let us calculate a family of Testing Curves for the set
$$S:= \{(x,y,z) : (x^2 + y^2 + z^2)x \geq z^2 \}. $$
From the form of $S$ it is easy to see that $\mult_S z^2 \geq 3$, hence $\mult_S$ is not the standard multiplicity at~$0$.

Let $\frkX =  \{ (x^2 + y^2 + z^2)x - z^2 = 0\}$. Then $\sigma=\beta_4\circ\beta_3\circ\beta_2\circ\beta_1$ is $\frkX$-admissible resolution where 
\begin{enumerate}
\item $\beta_1:(M_1,\frkX_1\cup H_1 , H_1)\to (\re^3,\frkX,0) $ is blow-up of $0$,
\item  $\beta_2:(M_2,\frkX_2\cup H_1\cup H_2 ,H_2)\to (M_1, \frkX_1\cup H_1 , p) $ is a blow-up of singular point $p$ on $\frkX_1\cap H_1$,
\item  $\beta_3:(M_3,\frkX_3 \cup H_1\cup H_2 \cup H_3 ,H_3)\to (M_2, \frkX_2\cup H_1\cup H_2 , C)$ is a blow-up of the smooth curve $C = H_1\cap \frkX_2$ where $\frkX_2$ is tangent to $H_1$
\item $\beta_4: (M_4,\frkX_4\cup H_1\cup\dots \cup H_4 ,H_4)\to (M_3, \frkX_3 \cup H_1\cup H_2 \cup H_3  ,  C') $ is blow-up of smooth curve $C'= H_1\cap H_3$ of points where $\frkX_3 \cup H_1\cup H_2 \cup H_3 $ is not SNC.
\end{enumerate}
We have $E(S)=\{i: S_\infty\cap H_i\neq \emptyset \} =\{2,3,4\}$ and all these components are $S$-essential. It suffices that we take testing curves parametrized by any open subsets of $H_i$ (in the sense of Lemma~\ref{lemma_normal_are_testing_curves}). 

We can choose local coordinates such that
\begin{align*}
\sgm_{|U_2}=\pi_2:U_2\to\re^3,\quad &  \pi_2(x,y,z) = ( x^2y, x, x^2z   ) \\
\sgm_{|U_3}=\pi_3:U_3\to\re^3,\quad &   \pi_3(x,y,z) = ( xz^2, xyz^2, x^2z^3   )\\
\sgm_{|U_4}=\pi_4:U_4\to\re^3,\quad &  \pi_4(x,y,z) = ( x^2 y^2z, x^2yz^2,x^3 y^2z^2  )
\end{align*}
where the sets $U_i$ are open, $H_i\cap U_i=\{x=0\}$ and $U_i$ intersects openly $H_i$ for $i=2,3,4$. In the local coordinates for $i=2,3,4$ we can assume $$U_i= (-\epsilon, \epsilon)\times V,\quad V= \left\{(y,z): |y|<1, |z|<1,  yz\neq 0\right\} .$$

Hence for parameters $(y,z)\in V  $ set
\begin{align*}
\gamma^2_{y,z}(t) =  ( yt^2, t, zt^2 ),  \quad
&  \gamma^3_{y,z}(t) =(z^2t, yz^2t , z^3t^2 ) ,\\
\gamma^4_{y,z}(t) &= ( y^2zt^2, yz^2t^2, y^2z^2t^3 )  
\end{align*}
for $t\in(0,\epsilon)$. Let $\Gamma_{y,z}$ be the curve with above three branches. Note that if $y<0$, then $\Gamma_{y,z}$  does not lie in $S$.  By construction, $(\Gamma_{y,z})_{(y,z)\in V}$ is the family of curves testing multiplicity relative to~$S$. We can for instance calculate
$$ \mult_S z = \min \left\{ 2,2, {3\over 2} \right\} = {3\over 2} \ \  {\it hence} \  \ \mult_Sz^2=3 . $$
\end{ex}

\subsection{Finite number of curves testing multiplicity of polynomials}
\label{sec_proof_of_reduction_thm}

\begin{lem}\label{lem_generic_with_parameter_is_generic}
Consider a function $\re^D\times U\ni (a,c) \to G(a,c)\in  \re$ which is polynomial in~$a$, $U\subset\re^m$. There exists finitely many parameters $c_1,\dots,c_k$, $k\leq D$, such that 
\begin{equation}\label{eqn_claim_for_G}
\forall_{a\in\re^D}\quad  \big( \   G(a,c_1)=\dots=G(a,c_k)=0\  \Leftrightarrow\ \forall_c G(a,c)= 0 \  \big)
\end{equation}
\end{lem}

\dowod From right to left implication is obvious. To prove the other direction consider algebraic sets $X_c=\{a\ :\ G(a,c)=0 \}$ and
$$X=\bigcap_{c\in U} X_c.$$
Due to Hilbert's basis theorem there exists finitely many $c_1,\dots, c_k$ such that $X=X_{c_1}\cap\dots\cap X_{c_k}$. To see that $c_1,\dots, c_k$ satisfy the claim let us assume to the contrary that~\eqref{eqn_claim_for_G} does not hold. Take $a_0$ such that $G(c_i, a_0)=0$ but $G(c,a_0)\neq 0 $ for some $c$. From the first part $a_0\in X_{c_1}\cap\dots\cap X_{c_k} = X$. Hence from definition of $X$ we have $G(a,c)=0$, which is a contradiction. 

The inequality $k\leq D$ follows from the fact that we can choose at most $D$ equations to define the variety in $D$-dimensional space.\koniecdowodu

Of course, Lemma~\ref{lem_generic_with_parameter_is_generic} can be written equivalently as: there exist finitely many $c_1,\dots, c_k$ such that
$$\forall_{a\in\re^D}  \left(  \exists_c G(a,c)\neq 0  \Leftrightarrow \exists_j G(a,c_j)\neq 0 \right).  $$

Now we can prove the main Theorem~\ref{thm_finite_mult} on reduction of dimension below. 
\begin{thm}\label{thm_finite_mult}
Let $S$ be a semianalytic subset of $\re^n$ fat at $0$. For any degree $d$ there exists a semialgebraic curve $\Gamma_d$ which tests multiplicity at $0$ with respect to $S$ for all polynomials from $\re_d[X]$ i.e.
$$ \big(\mult_S\big)_{|\re_d[X]} \ \equiv \ \big(\mult_{\Gamma_d}\big)_{|\re_d[X]} $$

Moreover, $\Gamma_d$ can be chosen so that its number of branches is not more than $d\cdot { {n+d}\choose{d}}\cdot N(S)$, where $N(S)$ depends on $S$. 
\end{thm}

\dowod 
By Theorem~\ref{thm_existence_rational_TC_mult} there exists a family $\{\Gamma_c \}_{c\in U}$ of  testing curves with $l$ branches  such that each analytically parametrized family of branches can be written as
 $$\arc_c = \sum\limits_{i=\ord\Gamma }^{\deg\Gamma} \overrightarrow{\varrho_i}(c) \cdot t^{i} ,\quad 0<t<<1.$$ 
If $U$ is finite, we get the claim. We will consider the case when $U$ is not finite.

Any polynomial $f=\sum_{|\alpha|\leq d} a_\alpha x^\alpha$ of degree $\leq d$ can be interpreted as the point $a\in\re^D$, where $D= { {n+d}\choose{d}}$ is the number of monomials of degree $\leq d$. 

Fix a polynomial $f$, then formally (i.e. without disregarding constant zero coefficients) we have
\begin{eqnarray} \label{eqn_coefficients_f_on_branch}
f(\arc_c(t))\ = \sum_{|\alpha|\leq d} a_\alpha  \left( \sum_i \varrho_{i,1}(c) t^{i}  \right)^{\alpha_1}\dots  \left( \sum_i \varrho_{i,n}(c) t^{i}  \right)^{\alpha_n} = \nonumber \\
 = \sum_{j=1}^{D'} \ \ \left( \sum_{  i | \alpha | = w_j } \beta_{i,\alpha}(c) a_\alpha   \right) \  t^{{w_j}},
\end{eqnarray}
where $w_1<\dots<w_{D'}$ and 
$$\{w_1,\dots,w_{D'}\} = \left\{\  i\cdot | \alpha | \ : \ |\alpha|\leq d, \ord\Gamma\leq i\leq\deg\Gamma \ \right\}.$$
Hence $D'\leq d\cdot(|\deg\Gamma-\ord\Gamma|+1) $. 

Consider the linear function on coefficients $G_c(a)=(G_1^c,\dots,G_{D'}^c):\re^D\to\re^{D'} $ with 
$$ G_j^c(a) = \sum_{ i | \alpha | =w_j } \beta_{i,\alpha}(c) a_\alpha. $$

Apply Lemma~\ref{lem_generic_with_parameter_is_generic} to $G_j$ so that we get finitely many parameters $c_1^j,\dots, c_{k_j}^j$ satisfying its claim. These parameters test vanishing of the coefficient $G_j^c(a)$ of $f\circ \arc_c$ in presentation~\eqref{eqn_coefficients_f_on_branch}. 
Without loss of generality, we can assume that 
$ \{c_1^i,\dots, c_k^i \} = \{c_1^j,\dots, c_k^j \} $ for every~$i,j\leq D'$ and denote this set as $ \{c_1,\dots,c_k\}$. We have $k\leq DD'$ and 
$$\ord_0 (f\circ\arc_c) = \min\{w_j: G_j^c(a)\neq 0 \} .$$

For generic $c_0$ we have 
$$ \ord_0 (f\circ\arc_{c_0}) = \inf_c \ord_0 (f\circ\arc_c) $$
and the minimum is attained, hence by  Lemma~\ref{lem_generic_with_parameter_is_generic} we get
$$  \inf_c \ord_0(f\circ\arc_c) = \min  \{w_j: \exists_{i} G_j^{c_i}(a)\neq 0 \} = \min \bigcup_{i}     \{w_j: G_j^{c_i}(a)\neq 0 \} = \mult_{\arc_{c_1}\cup\dots\cup\arc_{c_k} }f $$
Therefore for generic $c$ we have 
\begin{equation}\label{eqn_finitely_vs_all}
\ord_0 (f\circ\arc_c) = \mult_{\arc_{c_1}\cup\dots\cup\arc_{c_k} }f .
\end{equation}

Consider all $l$ branches of $\Gamma_c$ { and denote $c_1^j,\dots, c_{k_j}^j$ the parameters satisfying above equality~\eqref{eqn_finitely_vs_all} for consecutive families of branches $\arc_c^j$. 
Put $$\Gamma_d := \bigcup_{j=1,\dots, l}\bigcup_{i=1,\dots k_j} \arc_{c_i^j}^j.$$

}

The curve $ \Gamma_d $ tests the multiplicity relative to~$S$ of any polynomial $f$ of degree $\leq d$. 
{
Indeed,  for arbitrary polynomial $f$ we have 
$$\mult_{\Gamma_d} f = \min_i \mult_{\arc_{c_1}^i\cup\dots\cup\arc_{c_{k_i}}^i}f $$
and due to~\eqref{eqn_finitely_vs_all} we get
$$ \mult_{\Gamma_d}f = \min_i {\ord_0(f\circ\arc^i_c)\over \ord_0\|\arc^i_c\|} = \min_i {\ord_0(f\circ \arc^i_{c_0})\over \ord_0\|\arc^i_{c_0}\|}  $$
for generic $c$ and $c_0\notin \Xi_f$. Therefore, by definition of testing curves $\{\Gamma_c\}_{c\in U}$ we get 
$$\mult_{\Gamma_d}f=\mult_Sf $$
for any polynomial of degree $\leq d$.
}

Following the construction in the proof if the number of branches of $\{\Gamma_c\}$ is~$l$, then the number of branches at infinity of the curve $\Gamma_d$ is not more than $ldD(|\deg\Gamma-\ord\Gamma|+1)$. Recall that $l\leq E(S)$.
Hence the maximal number of branches of $\Gamma_d$ is not more than $dD\cdot N(S)$, where $D ={ {n+d}\choose{d}}$ and  $N(S)=l(|\deg\Gamma-\ord\Gamma|+1)$ depends only on~$S$ 
(up to choice of rectilinearization).    
\koniecdowodu

\begin{rk} 
In Theorem~\ref{thm_finite_mult} for fixed admissible rectilinearization of $\partial S$ one can choose $N(S)$ so that 
$$N(S)\ \leq \    (\deg\Gamma-\ord\Gamma+1) \ \cdot \ \#\{i \ :\ H_i\ {\rm is\ } S-{\rm essential} \},$$
compare Remarks~\ref{rk_number_of_branches_TC} and~\ref{rk_ord_deg_Gamma_TC}.
\end{rk}
In arbitrary dimension for a fixed set there may be several families of testing curves. We do not know whether for these families there exists a uniform minimal number of branches $l$ or minimal degree of truncation $\deg\Gamma$, since our method depends on the choice of resolution of singularities.

\begin{ex}
Let us continue Example~\ref{ex_3dim}. Fix degree $d$. We will show instances of~$\Gamma_d$, the curve testing multiplicity of all polynomials of degree $\leq d$. 

Note that $\deg\gamma_i^{y,z}-\ord\gamma_i^{y,z}=1$ for $i=2,3,4$ and we can assume $U_i = \re^2\setminus\{yz=0\}$. The curve~$\Gamma_d$ according to Theorem~\ref{thm_finite_mult} can be chosen as a union of branches $\gamma_{c_1}^i,\dots,\gamma_{c_{e(d)}}^i$ for $i=2,3,4$ with 
$$e(d)\leq 2 d  { {3+d}\choose{d}} = {d(d+1)(d+2)(d+3)\over 3}.$$
This estimation is not sharp. 
For instance, for $d=2$ it suffices to take $\Gamma_2$ as union of $12 $ branches $\gamma_{y_j,z_j}^i$, $i=2,3,4, j=1,2,3,4$ such that following determinants
 $$\det [1,y_j,z_j]_{j=1,2,3} ,\ \ 
 \det [y_j^2, z_j^2, y_jz_j]_{j=1,2,3}  ,\ \ 
 \det [1,y_jz_j, z_j, y_j^2z_j]_{j=1,\dots, 4} $$
are nonzero. For instance $(y_j, z_j )=(j, j+4)$, $j=1,2,3, 4$, satisfy this condition.
\end{ex}

\section{Stability  of relative multiplicity}\label{sec_stability}

In this section the main result is Theorem~\ref{thm_multiplicity_20} showing stability of multiplicity with respect to change of parameters in description of the set, provided the parameter is generic.  To this aim we will first prove with elementary methods that any sublevel set with respect to generic value of a mapping is a fat set, see Theorem~\ref{thm_fat_sublevel_at_p}, since we were unable to find a reference to such result in literature.

\subsection{Bifurcation values} 
Throughout this subsection let  $F:U\to \re^m$, with $U$ an open subset of a smooth manifold~$M$ of dimension $n$, be a $\mathcal{C}^\infty$-smooth mapping.

We say that a value $t\in\re^m$  is typical if $F$ is a $\mathcal{C}^\infty$-trivial fibration over~$t$ i.e. there exists  a neighborhood $\Upsilon$ of $t$ in $\re^m$ and a $\mathcal{C}^\infty$-diffeomorphism $\Phi: F^{-1}(t)\times \Upsilon \to F^{-1}(\Upsilon) $  such that $F(\Phi(x,s))=s$ for all $(x,s)\in F^{-1}(t)\times \Upsilon$. 
\label{bifurcation_values}
The set of values which are not typical is called the set ${\rm Bif }F$ of bifurcation values of $F$.

We use convention that the empty map is a trivial fibration. Therefore,  for $t\in\re^m\setminus \overline{F(U)}$  in a neighborhood~$\Upsilon$ of $t$ we have $F^{-1}(\Upsilon)=\emptyset$ and the fibration condition is met trivially. We consider such~$t$ to be also a typical value of $F$. Obviously $t\in \partial F(U) $ are not typical. 

Let us denote 
$${\rm Bif}_p F = \bigcap_{V:\ V {\rm\ is\ a\ nghbhd\ of\ }p} {\rm Bif}\left( F_{| V} \right),$$
the bifurcation values near a point $p\in M$. 
The set ${\rm Bif}_p F$ is closed in $\re^m$. 

In particular, if $M=\re^n$ and $p=\infty$, this is the standard definition of bifurcation values at infinity (where a neighborhood of infinity is a complement of a compact set). Obviously, ${\rm Bif}_p F\subset {\rm Bif}F$ for $p\in M$. If $p\in U$, then due to Ehresmann's Fibration Lemma we get ${\rm Bif}F_{|V} = \overline{\Crit F_{|V}}$ in some neighborhood $V$ of $p$, where $\Crit F_{|V} $ are the critical values of $F$ i.e. the set of images of points where rank of jacobian of $F_{|V}$ is less than $m$.

Although for smooth maps it may happen that ${\rm Bif}f=\re^m$, 
in the tame classes we are interested in, the following Bertini-Sard theorems hold.

\begin{fact}[\cite{KOS}, \cite{VERDIER}]\label{fact_Bertini_Sard}
Let $F:U\to \re^m$ be smooth, $U$ open subset of a smooth manifold $M$. 
\begin{itemize}
\item If $F$ is semialgebraic and $U=M=\re^n$, 
then ${\rm Bif}F = \Crit F \cup {\rm Bif}_\infty F$ is a closed subset of a semialgebraic set of codimension at least $1$.
\item If $F$ is proper real analytic and $M$ is a real analytic manifold, then  ${\rm Bif}F = \Crit F$ is a closed subset of a subanalytic set of codimension at least $1$.
\end{itemize}
\end{fact}

\subsection{Remarks on topology and fatness}
For any subset $A$ of a topological space denote by $\partial A = \overline{A}\setminus \Int(A) $ its border and by
$$ \NF(A) = A \setminus \overline{\Int(A) } $$
its set of nonfatness. Obviously $\NF(A)=\emptyset$ if and only if $A$ is fat.

\begin{lem}\label{lem_nonfat_pts_easy}
Let $A,B$ be sets in a metric space.
\begin{enumerate}
\item If $B\subset A$, then $$\NF(A)\subset \NF (B) \cup \left(\NF(A)\setminus B\right) .$$ 
 In particular, if $B$ is fat, then $\NF(A)\subset  \NF(A)\setminus B$
\item If $A$ is open and $B$ is fat, then $A\cap B$ is fat.
\end{enumerate}

\end{lem}
\dowod 
Let us prove (1). Indeed, let $x\in\NF(A)$. Let $x\notin \NF(B)$. If $x\in \overline{\Int (B) }$, then $x\in\overline{\Int (A) } $ contrary to the assumption. Hence $x\notin \overline{\Int (B) }\cup \NF(B)=B$ which proves the statement.

To prove (2) take any point $x\in A\cap B$.  From fatness of $B$ there exists a sequence $(x_n)\subset \Int B$ converging to $x$. From openess of $A$ it lies also in $A$. Hence $x\in \overline{\Int(A\cap B)}$. \koniecdowodu

\begin{prty}\label{prty_intersection_is_fat_induction_step}
Consider sets $A_0,\dots, A_m$ such that for 
any sequence $0\leq i_1<\dots<i_m\leq m$ of length $m$ the sets
$$A_{i_1}\cap\dots\cap A_{i_l}$$
are fat. Then
$$\NF(A_0\cap \dots \cap A_m) \  \subset \  \bigcap_{i=0}^m \partial A_i.$$
\end{prty}
\dowod
Let us remark that in general for any $A$ we have $\NF(A)=A\setminus \overline{\Int(A)}\subset \partial A$ and 
$\partial (A_0\cap\dots\cap A_m)\subset\bigcup_{i=0,\dots,m} \partial A_i$.

Assume $x\in \NF(A_0\cap\dots\cap A_m)$. Suppose that $x\notin \partial A_i$ for some $i$, say $i=0$.

Set  $A_1\cap\dots \cap A_m=A$ and $A_0=B$. By assumption $A$ is fat and $x\in \NF(A\cup B)\subset \partial A\cup \partial B$. If $x\in \partial A\setminus \partial B $, then from fatness of $A$ we get $x\in \overline{\Int (A)}$. There also exists an open set $V\ni x$ such that $V\cap \partial (B)=\emptyset$ and either $V\cap B=\emptyset$ or $V\subset \Int(B)$. In the first case $x\notin A\cap B$, a contradiction. In the second case, $x\in \Int(B)$ and $x\in \partial (A)$. Since $A$ is fat, there exists $(x_n)\subset \Int(A)$ such that $x_n\to x$. Since $x\in \Int(B)$, then $(x_n)\subset B$. Therefore $x_n\in \Int(A)\cap \Int (B) = \Int(A\cap B)$. Hence $x\in \overline{\Int(A\cap B)}$ and it is therefore a fat point of $A\cap B$ contrary to the assumption. Therefore $\NF(A\cap B)$ and $\partial (A)\setminus \partial (B)$ are disjoint. The arbitrary choice of $i$ ends the proof. 
\koniecdowodu

\subsection{Bifurcation values and fat sublevel sets} 
For  values $t,t'\in\re^m$ we denote by $[t,t']$ the closed hypercube $[t_1, t_1']\times\dots\times [t_m,t_m']$ 
and by $(t,t')$ its interior. Recall from  page~\pageref{page_order_definition_in_rectilinearization} that $t'\preceq t$ if and only if the set $[t',t]$ is nonempty. \label{page_notation_for_multi_segments}

\begin{lem}\label{lem_Ehresmann_to_loc_fat}
Let $F:U\to\re^m$ be smooth with $U$ open in manifold $M$. If $F$ is a trivial fibration over $\Upsilon\subset \re^m$, 
then for every $t,s\in \Upsilon$ we have 
$F^{-1}\left(\overline{(s,t)}\right) = \overline{ F^{-1}((s,t))}. $ 
\end{lem}
\dowod 
Whenever $F^{-1}\left(\overline{(s,t)}\right)$ is empty, the assertion holds trivially. So we may assume $m\leq n$ and $(s,t)\neq \emptyset$.  If $F$ is fibration over $\Upsilon$, there exists a diffeomorphism $\Phi$ of $F^{-1}(t)\times \Upsilon$ and $F^{-1}(\Upsilon)$. For any $s,t\in \Upsilon$ if $(s,t)\neq \emptyset$, we get
$$ \overline{F^{-1}((s,t))} =  \overline{ \Phi \left( F^{-1}(t)\times (s,t) \right) } = \Phi \left( F^{-1}(t)\times [s,t] \right) = F^{-1}([s,t])  $$
since $\Phi$ is open-closed. 
\koniecdowodu

By ${\rm proj}_{i_1,\dots,i_j}:\re^{m}\to \re^j$ denote the standard projection $x\to (x_{i_1},\dots,x_{i_j})$ onto $\re^j$. 

\begin{thm}\label{prop_fat_sublevel_at_p}\label{thm_fat_sublevel_at_p}
Let $F: U\to\re^m$ be a smooth map and $U$ open in a manifold~$M$. 
For $p\in M$ define
$$  B_pF=  \bigcup_{\substack{1\leq i_1<\dots<i_j\leq m\\ j=1,\dots, m}}  {\rm proj}_{i_1,\dots,i_j}^{-1} \left( {\rm Bif}_p({F_{i_1},\dots,F_{i_j}})\right). $$
For every $t\in\re^m\setminus B_pF$ the sublevel set
$$S_t := \{x\in U:\ F_1(x)\leq t_1,\dots, F_m(x)\leq t_m \}$$
is fat at $p$ i.e. germ of $S_t$ at $p$ is fat. 

Moreover, there exists a neighborhood $\Upsilon\subset \re^m$ of $t$ such that for all $s\in \Upsilon$ germ of $S_s$ at~$p$ is empty if and only if the germ $S_t$ at $p$ is empty.
\end{thm}
\dowod 
Note that $t\notin {\rm Bif}_p F$ if and only if $t$ is typical value of $F_{|V}$ for some open neighborhood $V$ of $p$. We will use induction with respect to the dimension $m$. 

Let $m=1$. Then  $B_pF={\rm Bif}_pF$. Take  a typical value $t$ of~$F_{|V}$. Hence for $s\preceq t$ such that $\|s-t\|$ is small enough apply Lemma~\ref{lem_Ehresmann_to_loc_fat} to get
\begin{eqnarray}
\overline{\Int(S_t\cap V)}\cap V\ =\  \overline{ F_{|V}^{-1}((-\infty,t))}\  =\  F_{|V}^{-1}((-\infty,s]) \cup \overline{   F_{|V}^{-1}((s,t)) } \ =\nonumber\\
 = \ F_{|V}^{-1}((-\infty,s])\cup  F_{|V}^{-1}([s,t])\  =\  S_t \cap V. \nonumber
\end{eqnarray}
Therefore,  the sublevel set $S_t$ is fat at $p$.

Consider $F=(F_0,\dots,F_m):U\to\re^{m+1}$. For any sequence $0\leq i_1<\dots<i_m\leq m$ by inductive assumption let  $B_p(F_{i_1},\dots,F_{i_m})$ satisfy the claim of the theorem for the mappings $(F_{i_1},\dots,F_{i_m})$. We have
$$B_pF = {\rm Bif}_pF \cup  \bigcup_{0\leq i_1<\dots<i_m\leq m}  {\rm proj}_{i_1,\dots,i_m}^{-1} \left( B_p(F_{i_1},\dots,F_{i_m})\right) $$
since projections commute.

Take $t\notin B_pF$. We will show that $S_t$ is fat at $p$. By definition of $B_pF$ there exists a neighborhood $V$ of~$p$ such that for every sequence $0\leq i_1<\dots<i_j\leq m$, with $j=1,\dots, m+1$, the value $(t_{i_1},\dots, t_{i_j})$ is a typical value of the mapping $(F_{i_1},\dots, F_{i_j})_{|V}$.

First, let us note that if $\max_{S_t\cap V}F_i\neq t_i$ for some $i$, say $i=0$, then $$S_t\cap V = \{F_0<t_0 \}\cap \{F_1\leq t_1,\dots, F_m\leq t_m\}\cap V.$$ Hence by the inductive assumption and Lemma~\ref{lem_nonfat_pts_easy}(3) we get the claim. Therefore, we may assume $\max_{S_t\cap V}F_i= t_i$ for all $i=0,\dots, m$.

Since both $\{ F_{0|V}\leq t_0 \}$ and $(F_{1|V},\dots, F_{m|V})^{-1}( (-\infty,t'] )$, with $t'=(t_1,\dots, t_m)$, are fat sets by assumption and case $m=1$ proved above, we get by Property~\ref{prty_intersection_is_fat_induction_step} and assumption $\max_{S_t\cap V}F_i= t_i$ that
$$ \NF(S_t \cap V)\subset \bigcap_{i=0}^m\{F_i=t_i\}\cap V .$$

But since $t\notin {\rm Bif}_pF$, then for $s\preceq t$ close enough to $t$ we have that the set $F_{|V}^{-1}\left(\overline{(s,t)}\right) $ is fat by  Lemma~\ref{lem_Ehresmann_to_loc_fat}. Since $F_{|V}^{-1}(t)\subset F_{|V}^{-1}\left(\overline{(s,t)}\right) \subset S_t\cap V $,  by Lemma~\ref{lem_nonfat_pts_easy}(2) we get that $NF(S_t\cap V)=\emptyset$. Hence the germ of $S_t$ is fat at $p$.  Induction ends the proof of the first part.

To prove the second claim, it suffices to note that $F$ is a trivial fibration over a neighborhood of any $t\notin B_pF$, because ${\rm Bif}_pF\subset B_pF$. In particular if the fiber $F^{-1}(t)$ is empty, then $F^{-1}(\Upsilon)$ is empty for a neighborhood $\Upsilon$ of $t$.
\koniecdowodu

\begin{cor}\label{cor_fat_sublevel}
Let $F: U\to\re^m$ be a smooth map and $U$ open in a manifold~$M$. 
Consider
$$  B F=  \bigcup_{\substack{1\leq i_1<\dots<i_j\leq m\\ j=1,\dots, m}}  {\rm proj}_{i_1,\dots,i_j}^{-1} \left( {\rm Bif}({F_{i_1},\dots,F_{i_j}})\right). $$
For every $t\in\re^m\setminus BF$ the sublevel set
$$S_t := \{x\in U:\ F_1(x)\leq t_1,\dots, F_m(x)\leq t_m \}$$
is fat.

Moreover, there exists a neighborhood $\Upsilon\subset \re^m$ of $t$ such that for all $s\in \Upsilon$ we have $S_s=\emptyset$  if and only if  $S_t=\emptyset$.
\end{cor}
\dowod Proof is the same as for Theorem~\ref{thm_fat_sublevel_at_p} but replacing ${\rm Bif}_pF$ by ${\rm Bif}F$.\koniecdowodu

\begin{cor}\label{cor_fat_is_generic}
Use notation and assumptions of Theorem~\ref{thm_fat_sublevel_at_p} and Corollary~\ref{cor_fat_sublevel}. Denote 
$$N_F:=\{ t\in \re^m:\ {\rm sublevel\ set\ } S_t {\rm \ is\ not\ fat\ in\ }M \}.$$
\begin{enumerate}
\item If $F$ is semialgebraic and $U=M=\re^n$, 
then $N_F$ is a semialgebraic set of codimension at most~$1$.
\item If $F$ is proper real analytic and $M$ is real analytic manifold, then  $N_F$ is a subanalytic set of codimension at most~$1$.
\end{enumerate}
If $U=M=\re^n$ put 
$$ N_F(\infty):=\{ t\in \re^m:\ {\rm sublevel\ set\ } S_t {\rm \ is\ not\ fat\ at\ }\infty \} .$$
\begin{itemize}
\item[(3)] If $F:\re^n\to\re^m$ is semialgebraic, then  $N_F(\infty)$ is a semialgebraic set of codimension at most~$1$.
\end{itemize}
\end{cor}
\dowod Note that $N_F$ and $N_F(\infty)$ are defined by  first order formulas. Moreover, $N_F\subset BF$ and $N_F(\infty)\subset B_\infty F$ by Theorem~\ref{thm_fat_sublevel_at_p} and Corollary~\ref{cor_fat_sublevel}. Use Fact~\ref{fact_Bertini_Sard}. \koniecdowodu

\subsection{Stability of relative multiplicity  under change of parameters} This section is devoted to the proof of the following

\begin{thm}[Stability of relative multiplicity]\label{thm_multiplicity_20}
Consider two real analytic mappings $f,g:\re^n\to\re^m$. For $t\in\re^m$ put
\begin{equation}\label{eqn_sublevel_set_mult}
S_t:=\{f_1+t_1g_1> 0,\dots, f_m+t_mg_m> 0  \}.
\end{equation}
There exists a nowhere dense subanalytic set $\mathcal{V}_{f,g}\subset \re^m$ such that 
$$\mult_{S_t}=\mult_{S_s}$$
provided $s,t$ lie in the same connected component of $\re^m\setminus \mathcal{V}_{f,g}$.
\end{thm}

Fix the mappings $f,g:\re^n\to\re^m$. Take $S_t $ of the form~\eqref{eqn_sublevel_set_mult} of statement of Theorem~\ref{thm_multiplicity_20}. Let 
$$ h:= \prod_{i,j,k,l\leq m}f_ig_j(f_k-g_l) $$
and $\frkX = \{h=0 \}$.
Obviously $\overline{S_0}\setminus S_0\subset \frkX$. Let $\sgm$ be the admissible rectilinearization of~$\frkX$ as in Section~\ref{section_resolution} with $E=H_1\cup\dots \cup H_r$, $F=H_{r+1}\cup \dots \cup H_p$. By choice of $h$, we can write $$\sgm^*f_i = I_{H_1}^{\mfE_1^i}\cdots I_{H_p}^{\mfE^i_p} {\quad \rm and\quad} \sgm^*g_j = I_{H_1}^{\mgE^j_1}\cdots I_{H_p}^{\mgE^j_p}, $$
where any $\mfE^i$ and $\mgE^j$ are comparable with respect to order $\preccurlyeq$.




Let us denote by $\mathcal{T}$ the set of all  values $t\in\re^m$ such that $(S_t)_{\infty}$ does not intersect~$E$ quasi-openly. Let us define the set $\mathcal{G}$ of good values as
$$ \mathcal{G} := \{  t\in \re^m\setminus \mathcal{T}\ :\    E(S_t)=E(S_s) {\ \rm for\ all\ }s{\rm \ in\  some\ nghbhd\ } U{\rm \ of \ }t  \} $$
(recall definition of $E(S)$ from page~\pageref{def_H_S}).
Recall that any $o$-minimal structure is closed under closure, interior, projection and boolean (first order) conditions. Therefore, the sets $\mathcal{T}$ and $\mathcal{G}$ are of respective category. 

\begin{lem}\label{lem_V}
The set  $\mathcal{G}$ contains an open dense subset of $\re^m$.
\end{lem}
\dowod 
Take a component $H$ of $E$.

Whenever $H$ is such that for every $j\leq m$ we have for multiplicities $\mfE^j\neq \mgE^j$, due to Border Lemma~\ref{lem_border_v3} there exists a neighborhood $V$ of $H$ such that the set $\sgm^{-1}(S_t\setminus\{h=0\})$ equals to a fixed union of connected components of $V\setminus (E\cup F)$ for all~$t$ in the same connected component of $\re^m\setminus\{\Pi_{i\in I} t_i=0 \} $ for some $I\subset \{1,\dots,r \}$. 
Moreover, $H\in E(S_t)$ if and only if $H\in E(S_s)$ for all $s,t $ in the same connected component of $\re^m\setminus\{\Pi_{i\in I} t_i=0 \} $. 
Therefore such component $H$ does not influence change between $E(S_t)$ and $E(S_s)$ as long as the parameters $t,s$ are in the same component of $\re^m\setminus \{\Pi_{i\in I} t_i=0 \}$.

Fix $H$ such that $\mfE^j = \mgE^j$ for $j\leq k$ and otherwise for $j>k$ (without loss of generality we changed the order of coordinates of $f$ and $g$). Consider 
$$\tilde{S}_t := \overline{\sgm^{-1}\left(\{f_1-t_1g_1\geq 0,\dots, f_k-t_kg_k \geq 0\}\setminus\{h=0\} \right)}.$$
Hence $(S_t)_\infty$ is fat in $H$ if  $(\tilde{S}_t)_\infty$ is fat in $H$.

Cover $H$ in $M$ by a family of open connected neighborhoods $U_x$ of $x\in H$ such that~$U_x$ is relatively compact and
$$ \sgm^*(f_i +t_i g_i )  = u_1^{\mfE_{j_1}^i}\cdots u_l^{\mfE_{j_l}^i} \Psi_{f_i} + t_i u_1^{\mgE_{j_1}^i}\cdots u_l^{\mgE_{j_l}^i} \Psi_{g_i} $$
on $U_x$, provided $x\in H_{j_1}\cap\dots\cap H_{j_l}$ for some components $H_{j_q}$ of $E\cup F$. By taking $U_x$ smaller we can assume  units $\Psi_{g_i}, \Psi_{f_i}$ are such that $\delta_x> |\Psi_{g_i}|, |\Psi_{f_i}|>\epsilon_x>0$ on~$U_x$ for some positive $\delta_x,\epsilon_x$.

If $U_x$ lies in the set of smooth points of $H$, 
possibly shrinking $U_x$, we get that the inequality $\sgm^*(f_i+t_ig_i)\geq 0$ on every component of $U_x\setminus H$ is given by $\Psi_{f_i}+t_i\Psi_{g_i} \trianglerighteq 0 $, where $\trianglerighteq \in \{\geq, \leq\}$ depends on the sign of $u^{\mgE}$ on the component.

Hence the set $(S_t)_\infty\cap U_x$ is a finite union of some sets of the form
$$ 
\bigcap_{j=1}^k \left\{  {\Psi_{f_j}\over \Psi_{g_j}} \trianglerighteq_j t_j  \right\} $$
where $(\trianglerighteq_j) \in\{\geq,\leq \}^k$ is some sequence of signs depending on parity of $\mgE$ and signs of $\Psi_{g_j}$. Recall that a finite union of fat sets is fat, hence without loss of generality we will assume equality of  $(S_t)_\infty\cap U_x$ with a set of such form.

For $j=1,\dots,k$ the functions $$F_j:= {\Psi_{f_j}\over \Psi_{g_j}} : H\cap U_x \to \re $$
are smooth and proper, defined on a smooth manifold $H\cap U_x$. Therefore from Corollary~\ref{cor_fat_is_generic} follows that for $F=(F_1,\dots, F_k)$ there exists a nowhere dense closed subset $N_{x,H}$ of $\re^k$ such that for all $t\notin N_{x,H}$ the set $(\tilde{S}_t)_\infty\cap H\cap U_x$ is fat in $H$.

Denote $$A_{x,H} := {\rm proj}^{-1}_{1,\dots,k}(N_{x,H})\subset\re^m.$$ The set $A_{x,H}$ is nowhere dense in $\re^m$, closed and for all $t\notin A_{x,H}$ the set $(S_t)_\infty\cap U_x$ is fat in $H_\reg$.

Let $x\in  H_{j_1}\cap\dots\cap H_{j_l}\cap H$ be a corner point of $H$ with  some components $H_{j_i}$ of $E\cup F$. Since $\mfE^i$ and $\mgE^i$ are comparable, we get 
\[   
\sgm^*(f_i +t_i g_i ) = 
     \begin{cases}
       u^{\mfE^i}   (\Psi_{f_i} + t_iu^{\mgE^i-\mfE^i}\Psi_{g_i}) & \mgE^i \succ \mfE^i \\
       u^{\mgE^i}   (u^{\mfE^i-\mgE^i}\Psi_{f_i} + t_i\Psi_{g_i}) &\mgE^i \prec \mfE^i \\
       u^{\mgE^i} (\Psi_{f_i} + t_i\Psi_{g_i})  &\mgE^i = \mfE^i\\
     \end{cases}
\]
where $u^{\mfE^i} = u_1^{\mfE_{j_1}^i}\cdots u_l^{\mfE_{j_l}^i} $. As before, if $\mgE^i \succ \mfE^i$, then (possibly shrinking $U_x$) on any connected component of $ U_x\setminus E$ the function $\sgm^*(f_i +t_i g_i ) $ does not vanish and its sign depends on the sign of $u^{\mfE^i} \Psi_{f_i} $. Analogously, if $\mgE^i \prec \mfE^i$, then the sign depends on the sign of $u^{\mgE^i}   t_i\Psi_{g_i}$. Therefore, if $\mgE^i \neq  \mfE^i$, then for all $t_i>0$ the set $ \overline{\sgm^{-1}\left(\{f_i-t_ig_i\geq 0\}\setminus\{h=0\} \right)}\cap U_x\cap E $ is the closure in $U_x$ of a fixed union of connected components of $U_x\setminus E$, same for $t_i<0$.

Consider all coordinates $f_i, g_i$ such that $\mfE^i=\mgE^i$, $i=1,\dots, k$ (after rearranging the numbering of the coordinates $f_i,g_i$) . Then the set $(S_t)_\infty \cap U_x$ is fat in $H$ if  $\tilde{S}_t \cap H\cap U_x$ is fat in $H$. As before, the functions $F_j:= {\Psi_{f_j}\over \Psi_{g_j}} : H\cap U_x \to \re $, $j=1,\dots, k$, are smooth and proper, hence there exists a closed nowhere dense set $A_{x,H}\subset \re^m$ (as in the case of smooth point, a lift of a set $N_{x,H}$) such that for $t\notin N_{x,H}$ the set $(S_t)_\infty \cap U_x$ is fat in $H$.

From the covering $(U_x)_{x\in H}$ of $H$  we can choose a finite subcovering $U_{x_1},\dots, U_{x_{N}}$. Then for all $t\notin \mathcal{V}_H:=  \{t_1\cdots t_m=0 \} \cup A_{x_1,H}\cup \dots \cup A_{x_N,H} $ the set $(S_t)_\infty\cap H$ is fat in $H$. Moreover, by the second part of Corollary~\ref{cor_fat_sublevel}, if we take $t\notin \mathcal{V}_H$ there exists a neighborhood $U\subset \re^m$ of $t$ such that $H\in E(S_t)$ if and only if $H\in E(S_s)$ for all~$s\in U$.


Put $\mathcal{V}:= \bigcup_{i=1}^r \mathcal{V}_{H_i} $, it is a nowhere dense closed set. Since every fat nonempty set contains a nonempty open set, then $\mathcal{T}\subset \mathcal{V} $. Moreover, if $t\notin \bigcup_{i=1}^r \mathcal{V}_i $, then there exists a neighborhood $U\subset \re^m$ of $t$ such that $H\in E(S_t)$ if and only if $H\in E(S_s)$ for all~$s\in U$. Therefore, the set $\mathcal{G}$ contains the set $\re^m\setminus  \mathcal{V} $ which proves the claim. \koniecdowodu

{\bf Proof of Theorem~\ref{thm_multiplicity_20}:}
We have $\overline{S_t\setminus \frkX} = \overline{S_t}$, because $S_t$ open and $\frkX$ nowhere dense.

Set $\mathcal{V}_{f,g} := \overline{ \re^m\setminus {\mathcal{G}} }$ of Lemma~\ref{lem_V}. It is nowhere dense and for any $t\notin \mathcal{V}_{f,g}$ and any component $H$ of $E$ the set $(S_t)_\infty\cap H$ is either empty or Zariski dense in $H$ i.e. $S_t$ intersects $E$ quasi-openly.

Since we can connect any two points $t,s$ in the same connected component of $\re^m\setminus \mathcal{V}_{f,g}$ by a compact curve, by standard argument and definition of $\mathcal{G}$ we get equality $E(S_t) = E(S_s)$ for any two point $t,s$ in the same connected component of $\re^m\setminus \mathcal{V}_{f,g}$.  
By Corollary~\ref{cor_equal_H_S} we get $\mult_{S_t}\equiv \mult_{S_s}$. \koniecdowodu

\section{Degree relative to a set}
In this section we  consider behavior of functions at infinity. For brevity 
 we restrict to the class of polynomials. More general results for polynomially bounded classes can be proved using direct methods of previous sections. 
 

\subsection{Degree with respect to a set}\label{section_degree}

 Let $S$ be a subset in $ \re^n$.  One characteristic of  behavior of a polynomial $f$ at infinity is the degree of $f$ relative to $S$ defined as

\begin{defi}[Rational degree relative to a set]\label{def_relative_degree}
 Let $S$ be unbounded. Put
$$ \rdeg_Sf:= \inf\left\{ {a\over b}: {|f(x)|^b \over \|x\|^a} {\rm \ is\ bounded\ on\  }S  {\rm \ outside \ some\ compact}  \right\} .$$
If $S$ is a bounded set, we assume $\rdeg_S\equiv -\infty$.
\end{defi}

Note that rational relative degree generalizes the relative degree of page~\pageref{def_relative_degree_integer}. Therefore,  for simplicity we will call the rational relative degree just the relative degree. Indeed,
\begin{rk}\label{rk_deg_and_rdeg}
We have $\deg_S =\max\left\{0,\lceil \rdeg_S \rceil\right\} $, where $\lceil\cdot\rceil$ is the ceiling function.
\end{rk}

One can alternatively write $$ \rdeg_Sf= \inf\left\{ {da\over b}: \exists_{g\in\re[X], d=\deg g} {|f(x)|^b \over |g(x)|^a} {\rm \ is\ bounded\ on\  }S  {\rm \ outside \ some\ compact }  \right\} $$
and when $S$ admits Curve Selection Lemma at infinity then
$$ \rdeg_Sf= \sup\left\{  {\deg (f\circ\arc)\over \deg\|\arc\|}  |\    \arc:[0,\infty)\to S, \lim_{t\to \infty}\|\arc\|=\infty  {\rm \ and\ }\arc{\rm \ is\ Laurent-Puiseux}  \right\}, $$
where by Laurent-Puiseux we mean a Laurent-Puiseux parametrization of an unbounded arc.

Let us remark that since polynomials are meromorphic at infinity, the (rational) relative degree function does not need to attain a minimum as in Example~\ref{ex_negative_deg} below.

\begin{ex}\label{ex_negative_deg}
Let $S=\{(x,y)\in\re^2 : |x^uy^w|\leq 1, y\geq 1 \}$ with $w,u\in\Na$. Then for any polynomial $f\in \re[X]$ independent of variable $Y$, we have
$$\rdeg_S f = -{w\over u}\ord f.$$
Indeed, it is easy to observe that on $S$ we have $|y| \cong \|(x,y)\|$ on $S$ and $|f| \lesssim |y|^{-{w\over u}\ord f} $. On the other hand, $\arc(t)=({1\over t^w} , t^u), t>>1,$ lies in $S$. We have $\deg (f\circ\arc)=-w\cdot \ord f$ and $\deg\arc=u$. Hence the assertion follows. 
\end{ex}

\subsection{Proofs of Theorems~\ref{main_thm_TC_deg} and~\ref{thm_degree_stability_20_integer} by inversion}\label{sec_inversion}

Note that the inversion in this subsection can be seen as a one-point compactification of punctured disc at infinity and imposing an appropriate metric, thus multiplicity and degree are essentially the same object. The degree of a polynomial, measure of its growth nearby $\infty$, can be read from the multiplicity at the point $\infty$. 
Therefore Theorems~\ref{main_thm_TC_deg} and~\ref{thm_degree_stability_20_integer} are a direct consequence of Theorems~\ref{thm_existence_testing_curves} and~\ref{thm_degree_stability_20} respectively, taking into account Remark~\ref{rk_deg_and_rdeg}.

Let us write $f\in\re[X]$ as a sum of homogeneous polynomials $f=h_d+\dots+h_0$, where $\deg f=d$, $\deg h_i = i$ and $h_i$ are forms. 

Denote the inversion $\iota(x):={x\over \|x\|^2}$ on $\re^n\setminus\{0\}$. Then
$$\iota^*f(x)  = {1\over\|x\|^{2d}}\left(h_d + \|x\|^2 h_{d-1}+\dots + \|x\|^{2d} h_0 \right).$$
Let us denote 
$$I(f):= \|x\|^{2\deg f}\iota^*f \ = \  h_d + \|x\|^2 h_{d-1}+\dots + \|x\|^{2d} h_0.$$

\begin{prop}[Duality of multiplicity and degree]\label{prop_duality_deg_mult}
Let $S\subset\re^n$ be semialgebraic, take polynomial $f\in\re[X]$ of degree $d$. We have
$$ \rdeg_S f = 2d-\mult_{\iota (S\setminus\{0\})} I(f).$$
\end{prop}
We leave proof of this Proposition to the reader. 

\begin{defi}[Testing Curves for relative degree]\label{def_testing_curves_deg}
Let $S$ be a set in $\re^n$. Let $(\Gamma_c)_{c\in U}$, $U\subset \re^m$, be a family of semialgebraic curves such that\begin{enumerate}
\item  there are $l$ connected components $\arc_1^c,\dots,\arc_l^c$ of $\Gamma_c$ for every $c$ and each component is parametrized as an arc
\item for each $i$ the family $(\arc_i^c)_c$ is parametrized analytically by $c$
\item for every polynomial $f:\re^n\to\re$ there exists a semialgebraic set $\Xi_f\subset U $ nowhere dense in $\re^m$ such that if $c\notin \Xi_f$, then $$\rdeg_S f = \max_{i} \deg (f\circ \arc_i^c ) $$ 
\end{enumerate}
We will call such a family $(\Gamma_c)$ the family of testing curves whenever $m<n$.
\end{defi}

Shortly, we call a family of curves $\{\Gamma_c\}$ testing curves for relative degree with respect to~$S$ if  for every polynomial~$f$ we have
$$\rdeg_S f=\deg_{\Gamma_c} f {\rm \ for \ generic \ }c$$
and the family is analytically parametrized by $c\in\re^m$, where $m<n$.

Note that genericity in the 3rd point of the definition has to depend on $f$. Precisely, for each $f$ the set of parameters $c\in U\subset \re^m$ such that the above equality does not hold is a globally analytic nowhere dense subset of $\re^{m}$ but its precise form depends on~$f$. On the other hand, for a given parameter $c$, the set of polynomials of degree~$d$ such that the above inequality does not hold is nowhere dense in $\re_d[X]$.

\begin{thm}[Existence of polynomial Testing Curves] 
For every semialgebraic set $S\subset \re^n$ fat  at  $\infty$ there exists an analytically parametrized family of semialgebraic testing curves for rational degree with respect to the set $S$.

Precisely, 
there exist integer numbers $q, \ord \Gamma,\deg \Gamma$ and testing curves of Definition~\ref{def_testing_curves_deg} with point (1) of the definition  replaced by
\begin{itemize}
\item[(1')] there are $l$ connected components $\arc_1^c,\dots,\arc_l^c$ of $\Gamma_c$ for every $c\in U$ such that each component is  parametrized by a finite Laurent-Puiseux series
$$ \arc_j^c(t) = \sum_{i=\ord\Gamma}^{\deg\Gamma} \overrightarrow{\varrho}_{j,i}(c) { t}^{i/q} ,$$
where $\overrightarrow{\varrho}_{j,i}:U\to\re^n$ are real analytic, $t>>1$.
\end{itemize}
\end{thm}
\dowod Use inversion of Section~\ref{sec_inversion} and apply Property~\ref{prop_duality_deg_mult} to Theorem~\ref{thm_existence_rational_TC_mult}. Normalize resulting testing curves, so that their degree is $1$.\koniecdowodu

\begin{thm}[Discrete family of Testing Curves]\label{thm_existence_testing_curves}
Let $S$ be a semialgebraic subset of~$\re^n$ fat at $0$. For any number $d$ there exists a semialgebraic curve $\Gamma_d$ with $l(d)$ branches at $\infty$ parametrized by Laurent-Puiseux finite series $\gamma_1,\dots, \gamma_{l(d)}$ 
such that
$$ \rdeg_S f =  \max_{i=1,\dots,l}  \deg \left(f\circ \arc_i \right)$$
for every $f\in\re[X]$, $\deg f=d$. 
Moreover,   
$$l(d)\leq d\cdot { {n+d}\choose{d}}\cdot N(S),$$ where $N(S)$ depends only on $S$. 
\end{thm}
\dowod Use inversion of Section~\ref{sec_inversion} and apply Property~\ref{prop_duality_deg_mult} to Theorem~\ref{thm_finite_mult}. Normalize parametrizations, so that their degree is $1$.
\koniecdowodu

Degree with respect to an arc $\Gamma$ can be seen as the degree of the univariate Laurent-Puiseux series $f\circ\Gamma$, where by abuse of notation $\Gamma$ is the parametrization of the arc. Calculation of $\deg ( f\circ\Gamma)$ is just symbolic.

\begin{thm}[Stability of degree for mappings]\label{thm_degree_stability_20}
Consider two polynomial mappings $f,g:\re^n\to\re^m$. For $t\in\re^m$ put
 $$S_t:=\{f_1+t_1g_1> 0,\dots, f_m+t_mg_m> 0  \}.$$ 

There exists a nowhere dense semialgebraic set $\mathcal{V}_{f,g}\subset \re^m$ such that 
$$\rdeg_{S_t}\equiv\rdeg_{S_s}$$
provided $s,t$ lie in the same connected component of $\re^m\setminus \mathcal{V}_{f,g}$.
\end{thm}
\dowod Use inversion of Section~\ref{sec_inversion}, apply Property~\ref{prop_duality_deg_mult} to Theorem~\ref{thm_multiplicity_20}.\koniecdowodu

\begin{ex}
Let us compute curves testing degree with respect to the set
$$ S:=\{ x\in\re^2  \ | \  -1\leq xy+y^2\leq 1 \} $$
fat at $\infty$. After inversion it can be written in the form
$$\iota(S)= \left\{  (xy+y^2)^2\leq (x^2+y^2)^4 \right\}$$
and is fat at $0$. Let $f=(xy+y^2)^2-(x^2+y^2)^4$ and $\frkX=\{ f=0\}$. One can check that after blow-up of zero the strict transform of $\frkX$ intersects the zero divisor at two points with multiplicity~$2$. Both points have to be blown-up twice so that the strict transform $\frkX$ intersects the zero divisor transversally. 

Denote $\beta (x,y) = (x,xy)$ and $l(x,y)=(x,y-1)$. Curves testing multiplicity with respect to $\iota(S)$ can be taken as $\beta\circ\beta\circ\beta (t,yt) = (t,t^3y)$ and $\beta\circ l\circ\beta\circ\beta (t,yt) = (t,t(t^2y-1))$ for $0<t<1$ and real parameter $y$.

Hence curves testing degree along $S$ can be chosen as union of branches
\begin{equation}\label{eqn_ex_2_dim_Tc}
 {(t,{y\over t}) \over 1+{y^2\over t^4}} \quad {\rm and}\quad  {(t,-t +{y\over t}) \over 1+{y^2\over t^4}} 
\end{equation}
for $t>>1$ and parameter $y\in \re$. 

Note that using methods of~\cite{MM_TC}, which are pertinent in 2-dimensional case, one can naturally reduce the number of branches by choosing testing curve to be curves given by one equation  $  xy+y^2 =c $ with $c\in\re$ (which is essentially the second branch in  presentation~\ref{eqn_ex_2_dim_Tc} when $y>>1$ and first branch when $|y|<<1$).
\end{ex}

\subsection{Filtration with respect to relative degree}
All properties and theorems of this section can be given also in terms of multiplicity, accordingly changing $\rdeg_S$ to $\mult_S$ and assumption that $S$ be fat at $\infty$ to fatness at $0$.

\begin{prty}\label{prty_only_rational_values_3_parts_deg}
Let $S$ be a semialgebraic subset of $\re^n$ fat at $\infty$.
\begin{enumerate} 
\item The relative degree function $\rdeg_S$
takes values on $\re[X]$ in a discrete set ${1/w}\Za$ with some $w=w(S)\in \Na$.
\item The set $\rdeg_S(\re_d[X])$
is finite for every $d$.
\end{enumerate}
\end{prty}
\dowod 
 Since the number of components of the zero divisor $E$ is finite, by Proposition~\ref{prop_how_to_compute_mult_via_divisor} there exists $w\in\Na$ such that $\rdeg_S f\in {1/w}\Za$ for all $f\in\re[X]$. Hence (1) holds.
 Proof of~(2) follows from rationality of testing curves, see Theorem~\ref{thm_existence_rational_TC_mult}.
\koniecdowodu

Let $$\mathcal{B}_q(S)=\{f\in\re[X]:\rdeg_S f\leq q \}$$ be the module of polynomials with relative degree not greater than $q$. Note that Property~\ref{prty_only_rational_values_3_parts_deg} point (3) does not mean that the module $\mathcal{B}_q(S)$ is of finite type. It can be seen in example $S= \{|xy|\leq 1 \}$, where $\mathcal{B}_0(S) = \re[X Y]$, hence $\mathcal{B}_0(S)$ contains polynomials of any degree. On the other hand, we can consider $\re_d[X]$ as the ${d+n}\choose{d}$ dimensional affine space and obtain the following.
\begin{prty}\label{prop_structure_of_degree_strata}
Let $S\subset \re^n$ be  a fat semialgebraic set. 
The set $\re_d[X]\cap \mathcal{B}_q(S)$ is a linear subspace of $\re_d[X]$ and 
$$\re_d[X] = \bigcup_{i=k}^K \{ f\in\re_d[X]\ :\ w\cdot\rdeg_S f=i  \}  $$
for some $K,w,k\in\Za$.
\end{prty}

Proof follows immediately from Property~\ref{prty_only_rational_values_3_parts_deg}. Note that in particular above shows that, unsurprisingly, for fixed~$S$ a  generic polynomial of degree~$d$ attains its degree on the set. 
Moreover, this gives a very strong quantifier elimination for the problem of relative multiplicities.

\begin{ex}
Let $S=\{0\leq xy-y\leq 1 \}$. Then testing curves for degree can be chosen as branches of $\gamma_c=\{xy-y=c \}$ at $\infty$. Let us describe relative degree grading on $\re_2[X,Y]= (a_0, a_{10}, a_{01}, a_{20}, a_{11}, a_{02})$. Put $D(d)=\{f\in\re_2[X,Y]: \rdeg_S f\leq d \}$. We have $$D(1) =\{a_{20}=a_{02}=0\},\ \ D(0)= D(1)\cap \{a_{10}=0, a_{01}+a_{11}=0 \},\ \  D(-1) = \{0 \}.$$
For instance $xy$ is of degree $1$ on $S$ and any power of $xy-y$ is of degree $0$ on $S$. Obviously, the grading does not depend on choice of testing curves.
\end{ex}

\subsection{Remarks}
Now let us present a basic closed version of Theorem~\ref{thm_degree_stability_20}.

\begin{thm}\label{thm_sublevel_set_stability}
Let $f=(f_1,\dots, f_m):\re^n\to\re^m$ be a polynomial mapping. Set
 $$S_t=\{f_1\geq t_1,\dots, f_m\geq t_m  \}.$$
There exists a nowhere dense semialgebraic set $\mathcal{V}_f\subset\re^m$ such that for any two $t,s$ in the same connected component of $\re^m\setminus\mathcal{V}_f$ we have
$$\rdeg_{S_t}\equiv \rdeg_{S_s} .$$
\end{thm}
\dowod 
Consider $p=\infty$ and take the set $B^\infty_f$ of Proposition~\ref{prop_fat_sublevel_at_p} and $\mathcal{V}_{f,1}$ for $Int(S_t)$ from Theorem~\ref{thm_degree_stability_20}. If $f:U\to\re^m$, $U\subset \re^n$, is polynomial, then ${\rm Bif }f$ is a subset of a proper algebraic set, see for instance~\cite{VERDIER}. Since ${\rm Bif}_\infty f\subset {\rm Bif}f$ and $\mathcal{V}_{f,1}$ is nowhere dense  due to Theorem~\ref{thm_degree_stability_20}, then $\mathcal{V}_f:= B^\infty_f \cup \mathcal{V}_{f,1}$ is a nowhere dense  set. Moreover,  for any $t,s$ in the same connected component of $\re^m\setminus \mathcal{V}_f$  there exists a compact $K$ such that $S_t\setminus K$ and $S_s \setminus K$ are fat. By Property~\ref{prty_degree_interior_vs_closure} and Theorem~\ref{thm_degree_stability_20} we get
$$ \rdeg_{S_t} \equiv  \rdeg_{Int S_t} \equiv \rdeg_{Int S_s} \equiv \rdeg_{S_s}. $$
This ends the proof.\koniecdowodu

Remark that the above Theorem~\ref{thm_sublevel_set_stability} generalizes a result of~\cite{KMS}. But in~\cite{KMS} for principal semialgebraic subsets of the real plane, i.e. $n=2 $ and $m=1$, it was additionally shown  that $$\mathcal{V}_f\subset {\rm Bif}_\infty f.$$
It would be interesting to establish this relation for $n>2$. Checking if $\rdeg_{S_t} = \rdeg_{S_s}$ involves only countably many polynomials in one variable due to reduction of dimension in Theorem~\ref{thm_existence_testing_curves},  or countably many  linear equations, see Proposition~\ref{prop_structure_of_degree_strata}, hence it is relatively easy. 
If indeed $\mathcal{V}_f\subset {\rm Bif}_\infty f$ for $n>2$ was true, then any generalized critical value at which the relative degree function changes would be necessarily a bifurcation value at infinity. 

Note that in Example~\ref{ex_unstab_neq_bifurcation} below we show that unfortunately the instability values $\mathcal{V}_f$ can be a proper subset of bifurcation values at infinity.

\begin{ex}
\label{ex_unstab_neq_bifurcation}
Let $\mathcal{V}_f$ be the smallest set such that claim of Theorem~\ref{thm_sublevel_set_stability} holds. Let $B=x(xy-1)$ be the Broughton polynomial, $a>0$ and  $$f(x,y):= (B(x,y)-a)^2 .$$
Put $S_t=\{f\geq t \}$. Then $\mathcal{V}_f =\{0\} \varsubsetneq {\rm Bif}_\infty f =\{0,a^2 \}$ as the link of $S_t$ at infinity changes at $t=a^2$. To show the statement,  one can use explicit methods of~\cite{MM_TC} to show testing curves for $S_t$, $t>0$, are  $\{x^2y=c, y\geq 1 \}_{c\in\re}$ and $\{x^2y-x=c, x\geq 1 \}_{c\in\re}$. They do not depend on $t>0$, in particular relative degree does not change near $a^2$.
\end{ex}

\begin{prop}\label{prop_stability_closed_set_one_f_for_mult}
Fix two polynomials $f,g$ and set  $S^\geq_t:= \{ f+tg\geq 0 \}$, $t\in\re$. Then the family of functions $\{\rdeg_{S_t}\}_{t\in\re}$ is finite.
\end{prop}
\dowod Using inversion it is sufficient to prove the following:

\em Let $S^\geq_t:= \{ f+tg\geq 0 \}$, $t\in\re$. Then the family of functions $\{\mult_{S_t}\}_{t\in\re}$ is finite.\em

As before, denote $\frkX=\{fg(f-g)=0 \}$. 
Fix $\mathcal{G}$ as in Lemma~\ref{lem_V} for $S_t$. Note that $(S_t)_\infty = (S_t^\geq)_\infty$. Hence  $S_t^\geq$ intersects $E$ quasi-openly if and only if $S_t$ intersects~$E$ quasi-openly. Moreover, for every $t,s$ in the same connected component of $\mathcal{G}$ we have $E(S_t^\geq)= E(S_s^\geq)$.

Note that 
\[   
S^\geq_t\cap \frkX = 
     \begin{cases}
       \{f=0, g\leq 0 \} \cup \{ f\geq 0, g=0 \}\cup \{f=g,f\leq 0 \} &{\rm for\ }t<0 \\
       \{f=0, g\geq 0 \} \cup \{ f\geq 0, g=0 \}\cup \{f=g,f\leq 0 \} &{\rm for\ }t\in (0,1) \\
       \{f=0, g\geq 0 \} \cup \{ f\geq 0, g=0 \}\cup \{f=g,f\geq 0 \} &{\rm for\ }t>1 \\
     \end{cases}
\]
Take the finite set $\mathcal{V}:= \{0,1\}\cup (\re\setminus \mathcal{G})  = \{t_1<\dots<t_N \}.$ Set $t_0=-\infty, t_{N+1}=+\infty$. By Corollary~\ref{cor_equal_H_S} for every $t,s\in (t_i, t_{i+1}) $ we have 
$$\mult_{S_t^\geq} = \mult_{S_s^\geq}.$$
Hence there are at most $2\#\mathcal{V}+1$ functions that are equal to a multiplicity relative to $S_t^\geq$.

(Alternatively, one can replace this simple proof by following proof of Theorem~\ref{thm_testing_curves_mult_v2} and showing that $(S^{\geq}_t)_F$ is either constant for $t>0$, or $\{i: (S^{\geq}_t)_F\cap H_i\neq\emptyset \} = E(S^{\geq}_t)$).\koniecdowodu

Of course, we cannot expect the set $S_t^\geq$ of Proposition~\ref{prop_stability_closed_set_one_f_for_mult} above to be fat at $0$ for generic~$t\in\re$. Indeed, for instance if $f=hf', g=hg'$ with $-h,f',g' \geq 0$ on $\re^n$, then for any $t>0$ we have $\{f+tg\geq 0 \} = \{h=0 \}$. 
Let us illustrate this with an example.
\begin{ex}
Let $S_t=\{ x^2-y^2\leq t (xy)^2 \}$. For $t<0 $ relative degree $\rdeg_{S_t}$ is the degree with respect to variable $Y$  whereas   $\rdeg_{S_t}$ is the standard degree for $t\geq0$.
\end{ex}

Note that from Theorem~\ref{thm_degree_stability_20} we get immediately
\begin{cor}\label{cor_rdeg_is_generically_a_finite_family}
Under assumptions and notations of Theorem~\ref{thm_degree_stability_20}, there exists an open dense set $\mathcal{G}$  of the parameter space $\re^m$ and a finite collection of functions $\rdeg_i:\re[X]\to\mathbb{Q}$, $i=1\dots,q$, such that for every $t\in \mathcal{G}$ exists $i$ such that 
$$ \rdeg_{S_t}\equiv\rdeg_i. $$
\end{cor}

This corollary raises a question whether there exists a stratification of the whole parameter space such that multiplicity is constant on strata. Let us look at following degenerate example. 
\begin{ex}\label{ex_there_can_be_moduli}
Consider $S_{t_1,t_2} = \{f\geq t_1, f\leq t_2 \} $. On the line $\{t_1=t_2\}\subset \re^2$ we have
$$ \forall_{t\neq s}\quad \mult_{S_{t,t}}\neq \mult_{S_{s,s}}.$$
Indeed, $\mult_{S_{t,t}}(f-t)=\infty\neq \mult_{S_{s,s}}(f-t)$.
\end{ex}

Hence in case of sets described by  two and more inequalities, there may be moduli of relative degree.

\section*{Acknowledgements}
V. Grandjean was partially supported by FUNCAP/CAPES/CNPq grant 305614/2015 and grant 306119/2018-8; M. Michalska was partially supported by NCN grant 2013/09/D/ST1/03701 and FAPESP grant 2015/08149-9.

\bibliographystyle{alpha}
\bibliography{Grandjean_Michalska_Degree_and_Multiplicity}

\end{document}